\documentclass[10pt]{article}
\usepackage[dvips]{graphicx}
\usepackage{graphicx}
\usepackage{amsmath}
\usepackage{amssymb}
\usepackage{amsthm}
\usepackage{amsfonts}
\usepackage{latexsym}
\usepackage{url}
\theoremstyle{plain}
\newtheorem{theorem}{Theorem}[section]

\newtheorem{lemma}[theorem]{Lemma}
\newtheorem{corollary}[theorem]{Corollary}
\newtheorem{proposition}[theorem]{Proposition}

\theoremstyle{definition}
\newtheorem{remark}{Remark}[section]
\newtheorem{example}{Example}[section]
\newtheorem*{notation}{Notation}

\theoremstyle{remark}
\newtheorem*{note}{Note}

\newcommand{\cB}{{\mathcal B}}

\newcommand{\cS}{{\mathcal S}}
\newcommand{\cA}{{\mathcal A}}

\newcommand{\cD}{{\mathcal D}}
\newcommand{\cP}{{\mathcal P}}

\newcommand{\cC}{{\mathcal C}}

\newcommand{\QQ}{{\mathbb Q}}
\newcommand{\NN}{{\mathbb N}}
\newcommand{\PP}{{\mathbb P}}

\newcommand{\cAstar}{\mathcal{A}^*}
\newcommand{\cAplus}{\mathcal{A}^+}
\newcommand{\empt}{\varepsilon}
\newcommand{\cAw}{\mathcal{A}^{\omega}}
\newcommand{\cAinf}{\mathcal{A}^{\infty}}

\newcommand{\rev}{\widetilde}

\newcommand{\vqed}{\vspace{-1cm} \flushright{\qedsymbol}}

\newcommand{\pref}{\prec_p}
\newcommand{\suff}{\prec_s}

\newcommand{\ov}{\overline}
\newcommand{\bs}{\mathbf{s}}
\newcommand{\bx}{\mathbf{x}}
\newcommand{\bt}{\mathbf{t}}

\newcommand{\bz}{\mathbf{z}}

\newcommand{\by}{\mathbf{y}}

\numberwithin{equation}{section}

\author{Amy Glen\footnotemark[1]}
\title{Powers in a class of $\cA$-strict standard episturmian words}
\date{August 27, 2006} 
\begin{document}
\normalsize
\normalsize \maketitle \footnotetext[1]{E-mail:
\texttt{amy.glen@gmail.com}}
\begin{center}
\vspace{-0.7cm} School of Mathematical Sciences,
Discipline of Pure Mathematics, The University of Adelaide, \\
South Australia, Australia, 5005
\end{center}
\vspace{0.5cm} \hrule 
\begin{abstract} 
This paper concerns a specific class of {\em strict standard episturmian words} whose {\em directive words} resemble those of {\em characteristic Sturmian words}. In 
particular, we explicitly determine all integer powers occurring in such infinite words, extending recent results of Damanik and Lenz (2003), who studied powers in Sturmian words. The key tools in our analysis are canonical decompositions and a generalization of \emph{singular words}, 
which were originally defined for the ubiquitous \emph{Fibonacci word}. Our main 
results are demonstrated via some examples, including the 
\emph{$k$-bonacci word}, a generalization of the Fibonacci word to a 
$k$-letter alphabet ($k\geq2$). 
\vspace{0.1cm}   
\vspace{0.1cm} 
\end{abstract}
\hrule

\section{Introduction}

Introduced by Droubay, Justin and Pirillo \cite{xDjJgP01epis}, 
\emph{episturmian words} are an interesting natural generalization of the well-known family of
\emph{Sturmian words} 
(aperiodic infinite words of minimal complexity) to an 
arbitrary finite alphabet. Episturmian words share many properties with Sturmian words and include the well-known \emph{Arnoux-Rauzy sequences}, the study of which began in  \cite{pAgR91repr} (also see \cite{ jJgP02onac, rRlZ00agen} for example). 

In this paper, the study of episturmian words is continued in more detail. In particular, for a specific class 
of episturmian words (a typical element of which we shall  
denote by $\bs$), we will explicitly determine all the integer powers 
occurring in its constituents.  
This has recently been done in \cite{dDdL03powe} for Sturmian words, which are exactly the aperiodic episturmian words over a 2-letter alphabet. 

A finite word $w$ has a (non-trivial) integer \emph{power} 
in an infinite word $\bx$ if $w^p = ww\cdots w$ ($p$ times) is a 
\emph{factor} of $\bx$ for some integer $p\geq2$. Here, 
our analysis of powers occurring in episturmian words $\bs$ hinges on 
canonical decompositions in terms of their `building blocks'. 
Another key tool is a generalization of \emph{singular words}, which 
were first defined in \cite{zWzW94some} 
for the ubiquitous \emph{Fibonacci word}, and later extended to Sturmian words 
in \cite{gM99lynd} and the \emph{Tribonacci sequence} in \cite{bTzW03some}. 
Our generalized singular words will prove to be useful in the study of 
factors of episturmian words, just as they have been for Sturmian words.

This paper is organized as follows. After some preliminaries (Section 
\ref{S:prelim}), we define, in Section \ref{S:family}, a restricted class of 
episturmian words upon which we will focus for the rest of the 
paper. A typical element of this class will be denoted by $\bs$. 
In Section \ref{S:g_singular}, we prove some simple results, which 
lead us to a generalization of \emph{singular words} for
episturmian words $\bs$. The \emph{index}, i.e., 
maximal \emph{fractional power}, of the building blocks of $\bs$ is 
then studied in Section \ref{S:index}. 
Finally, in Section \ref{S:powers_in_s}, we determine all squares 
(and subsequently higher powers) occurring in $\bs$. Our main results are 
demonstrated via some examples, including the \emph{$k$-bonacci word}, a 
generalization of the Fibonacci word to a $k$-letter alphabet ($k\geq2$).

\section{Definitions and notations} \label{S:prelim}

\subsection{Words}

Let $\cA$ denote a finite alphabet. A (finite) \emph{word} is an
element of the \emph{free monoid} $\cAstar$ generated by $\cA$, in
the sense of concatenation. The identity $\empt$ of $\cAstar$ is
called the \emph{empty word}, and the \emph{free semigroup},
denoted by $\cAplus$, is defined by $\cAplus :=
\cAstar\setminus\{\empt\}$. An \emph{infinite word} (or simply 
\emph{sequence}) $\bx$ is a sequence indexed by $\NN$ with values in $\cA$, 
i.e., $\bx = x_0x_1x_2\cdots$, where each $x_i \in \cA$. The set of 
all infinite words over $\cA$ is denoted by $\cAw$, and we define 
$\cAinf := \cAstar\cup\cAw$. 
If $u$ is a non-empty finite word, then $u^\omega$ denotes the \emph{purely
periodic} infinite word $uuu\cdots$.

If $w = x_{1}x_{2}\cdots x_{m} \in \cAplus$, each 
$x_{i} \in \cA$, the \emph{length} of $w$ is $|w| = m$ and we denote by 
$|w|_a$ the number of occurrences of a letter $a$ in $w$. (Note that 
$|\empt| = 0$.) The \emph{reversal} of $w$  
is $\rev{w} = x_{m}x_{m-1}\cdots x_{1}$, and if 
$w = \rev{w}$, then $w$ is called a \emph{palindrome}.

A finite word $w$ is a \emph{factor} of $z \in \cAinf$ if $z =
uwv$ for some $u \in \cAstar$, $v \in \cAinf$, and we write $w
\prec z$. Further, $w$ is called a \emph{prefix} (resp.~\emph{suffix}) of $z$ if $u = \empt$ (resp.~$v = \empt$), and we
write $w \pref z$ (resp.~$w \suff z$). 

An infinite word $\bx \in
\cAw$ is called a \emph{suffix} of $\bz \in \cAw$ if there exists a
word $w \in \cA^+$ such that $\bz = w\bx$. A factor $w$ of a
word $z \in \cAinf$ is \emph{right} (resp.~\emph{left})
\emph{special} if $wa$, $wb$ (resp.~$aw$, $bw$) are factors of $z$
for some letters $a$, $b \in \cA$, $a \ne b$.

For any word $w \in \cAinf$, $\Omega(w)$ denotes the set of all its
factors, and $\Omega_n(w)$ denotes the set of all factors of
$w$ of length $n \in \NN$, i.e., $\Omega_n(w) := \Omega(w)
\cap \cA^n$ (where $|w|\geq n$ for $w$ finite). Moreover, the \emph{alphabet} of $w$ is Alph$(w)
:= \Omega(w) \cap \cA$ and, if $w$ is infinite, we denote by Ult$(w)$ the set of
all letters occurring infinitely often in $w$.  Two infinite words $\bx$, $\by \in \cAw$ are said to be \emph{equivalent} if $\Omega(\by) = \Omega(\bx)$, i.e., if $\bx$ and $\by$ have the same set of factors.

Let $w = x_1x_2\cdots x_m \in \cAstar$, each $x_i \in \cA$, and
let $j \in \NN$ with $0 \leq j \leq m - 1$. The \emph{$j$-th conjugate} of 
$w$ is the word $C_j(w) := x_{j+1}x_{j+2}\cdots x_{m}x_1x_2\cdots x_j$, and 
we denote by $\cC(w)$ the conjugacy class of $w$, i.e., $\cC(w) :=
\{C_j(w) ~: ~0 \leq j \leq |w| - 1 \}$.  
Observe that if $w$ is \emph{primitive} (i.e., not a power of a shorter
word), then $w$ has exactly $|w|$ distinct conjugates.

The \emph{inverse} of $w \in \cAstar$, written $w^{-1}$, is
defined by $ww^{-1} = w^{-1}w = \empt$. It must be emphasized that
this is merely formal notation, i.e., for $u, v, w \in \cAstar$, the
words $u^{-1}w$ and $wv^{-1}$ are defined only if $u$ (resp.~$v$)
is a prefix (resp.~suffix) of $w$.

A \emph{morphism on} $\cA$ is a map $\psi: \cAstar \rightarrow
\cAstar$ such that $\psi(uv) = \psi(u)\psi(v)$ for all $u, v \in
\cAstar$.  It is uniquely determined by its image on the alphabet
$\cA$.

\subsection{Episturmian words}

An infinite word $\bt \in \cAw$ is \emph{episturmian} if $\Omega(\bt)$ is closed under
reversal and $\bt$ has at most one right (or equivalently left) special factor of each length.  
Moreover, an episturmian word is \emph{standard} if all of its left special factors are prefixes of
it. 

Standard episturmian words are characterized in \cite{xDjJgP01epis} using the concept of the \emph{palindromic right-closure} $w^{(+)}$ of a finite word $w$, which is the (unique) shortest palindrome having $w$ as a prefix (see \cite{aD97stur}).    
Specifically, an infinite word $\bt \in \cAw$ is standard episturmian if and only if  there exists an infinite word $\Delta(\bt) = x_1x_2x_3\ldots$, each $x_i
\in \cA$, called the \emph{directive word} of $\bt$, such that the infinite sequence of palindromic prefixes $u_1 =
\empt$, $u_2$, $u_3$, $\ldots$ of $\bt$ (which exists by results in \cite{xDjJgP01epis}) is given by 
\begin{equation} \label{eq:02.09.04}
  u_{n+1} = (u_nx_n)^{(+)}, \quad 
  n \in \NN^+.
\end{equation}
\begin{note} For any $w \in \cA^+$, $w^{(+)} = wv^{-1}\rev{w}$ where $v$ is the longest palindromic suffix of $w$.
\end{note}

An important point is that a standard episturmian word 
$\bt$ can be constructed as a limit of an infinite sequence of its palindromic prefixes, i.e., 
$\bt = \lim_{n\rightarrow\infty} u_n$. 

Let $a \in \cA$ and denote by $\Psi_a$ the morphism on $\cA$ defined by 
\[
\Psi_a : \begin{array}{lll}
                a &\mapsto &a \\
                x &\mapsto &ax \quad \mbox{for all $x \in \cA\setminus\{a\}$}.
                \end{array} 
\]
Another useful characterization of standard episturmian words is the following (see \cite{jJgP02epis}). An infinite word $\bt \in \cAw$ is standard episturmian with directive word $\Delta(\bt) = x_1x_2x_3\cdots$ ($x_i \in \cA$) if and only if there exists an infinite sequence of infinite words
$\bt^{(0)} = \bt$, $\bt^{(1)}$, $\bt^{(2)}$, $\ldots$ such that $\bt^{(i-1)} = \Psi_{x_i}(\bt^{(i)})$ for all $i \in \NN^+$. Moreover, each $\bt^{(i)}$ is a standard episturmian word with directive word 
$\Delta(\bt^{(i)}) = x_{i+1}x_{i+2}x_{i+3}\cdots$, the \emph{$i$-th shift} of $\Delta(\bt)$. 

To the prefixes of the directive word $\Delta(\bt) = x_1x_2\cdots$, we associate the morphisms 
\[
  \mu_0 := \mbox{Id}, \quad \mu_n := \Psi_{x_1}\Psi_{x_2}\cdots\Psi_{x_n}, \quad n \in \NN^+, 
\]
and define the words 
\[
  h_n := \mu_n(x_{n+1}), \quad n \in \NN,
\]
which are clearly prefixes of $\bt$. We have the following useful formula \cite{jJgP02epis}
\[
  u_{n+1} = h_{n-1}u_{n};
\]
and whence, for $n > 1$ and $0 < p < n$,
\begin{equation} \label{eq:u_n}
  u_n = h_{n-2}h_{n-3}\cdots h_1h_0 = h_{n-2}h_{n-3}\cdots h_{p-1}u_p.
\end{equation}

Some useful properties of the words $h_n$ and $u_n$ are given by 
the following lemma.

\begin{lemma} \emph{\cite{jJgP02epis}} \label{L:h_n&u_n}
For all $n \in \NN$,
\begin{itemize}
\item[\emph{(i)}] $h_n$ is a primitive word$;$
\item[{\em (ii)}] $h_n = h_{n-1}$ if and only if $x_{n+1} = x_n;$
\item[{\em (iii)}] if $x_{n+1} \ne x_n$, then $u_n$ is a proper prefix
of $h_n$.
\end{itemize} \vqed
\end{lemma}

Two functions can be defined with regard to positions of letters in a 
given directive word.  
For $n \in \NN^+$, let $P(n) = \sup\{p < n ~: ~x_p = x_n\}$ 
if this integer exists, $P(n)$ undefined otherwise. Also, let $S(n) =
\inf\{p > n ~: ~x_p = x_n \}$ if this integer exists, $S(n)$
undefined otherwise. By the definitions of palindromic closure and
the words $u_n$, it follows that $u_{n+1} = u_nx_nu_n$ (whence
$h_{n-1} = u_nx_n$) if $x_n$ does not occur in $u_n$, and $u_{n+1}
= u_nu_{P(n)}^{-1}u_n$ (whence $h_{n-1}u_{P(n)}= u_n$) if $x_n$
occurs in $u_n$. Thus, if $P(n)$ exists, then
\begin{equation} \label{eq:h_n-1}
  h_{n-1} = h_{n-2}h_{n-3}\cdots h_{P(n)-1}, \quad n \geq 1.
\end{equation}

\subsubsection{Strict episturmian words}

A standard episturmian word $\bt\in \cA^\omega$, or any equivalent (episturmian) word,   
is said to be \emph{$\cB$-strict} (or $k$-\emph{strict} if $|\cB|=k$, or {\em strict} if $\cB$ is understood) if 
Alph$(\Delta(\bt)) =$ Ult$(\Delta(\bt)) = \cB \subseteq \cA$.  In particular, a standard episturmian word 
over $\cA$ is $\cA$-strict if every letter in $\cA$ occurs infinitely many times in its directive word.  
The $k$-strict episturmian words have \emph{complexity} $(k - 1)n + 1$ for 
each $n \in \NN$ (i.e., $(k-1)n+1$ distinct factors of length $n$ for each $n\in \NN$). Such words are exactly the $k$-letter {\em Arnoux-Rauzy sequences}.

\subsubsection{Return words} \label{SS:return}

Let $\bx \in \cAw$ be \emph{recurrent}, i.e., any factor 
$w$ of $\bx$ occurs infinitely often in $\bx$. A  
\emph{return word}  of $w \in \Omega(\bx)$ is a factor of 
$\bx$ that begins at an occurrence of $w$ in $\bx$ and ends exactly before the next occurrence of $w$ 
in $\bx$. Thus, a return word of $w$ is a non-empty factor $u$ of $\bx$ such that $w$ is a prefix of 
$uw$ and $uw$ contains two distinct occurrences of $w$. This notion was introduced independently by Durand \cite{fD98acha}, and Holton and Zamboni \cite{cHlZ99desc}. 

Episturmian words are recurrent and, according to 
\cite[Corollary 4.5]{jJlV00retu}, each factor of an $\cA$-strict episturmian 
word has exactly $|\cA|$ return words.

\section{A class of strict standard episturmian words} \label{S:family}

Given any infinite sequence $\Delta = x_1x_2x_3\cdots$ over a
finite alphabet $\cA$, we can define a standard episturmian word
having $\Delta$ as its directive word (using \eqref{eq:02.09.04}).
In this paper, however, we shall only consider a specific family of
$\cA$-strict standard episturmian words.

Let $\cA_k$ denote a $k$-letter alphabet, say $\cA_k =
\{a_1,a_2,\ldots, a_k\}$, and suppose $\bt$ is a 
standard episturmian word over $\cA_k$. Then the directive word of
$\bt$ can be expressed as:
\[
  \Delta(\bt) = a_1^{d_1}a_{2}^{d_2}\cdots a_k^{d_k}a_{1}^{d_{k+1}}
           a_{2}^{d_{k+2}}\cdots a_{k}^{d_{2k}}a_{1}^{d_{2k+1}}\cdots,
\]
where the $d_i$ are non-negative integers. 
In what follows, we restrict our 
attention to the case when all $d_i > 0$; that is, we shall only study the class of $k$-strict standard episturmian
words $\bs \in \cA_k^\omega$ with directive words of the
form:
\begin{equation} \label{eq:dir_seq}
  \Delta = a_1^{d_1}a_{2}^{d_2}\cdots a_k^{d_k}a_{1}^{d_{k+1}}
           a_{2}^{d_{k+2}}\cdots a_{k}^{d_{2k}}a_{1}^{d_{2k+1}}\cdots,
         \quad d_i > 0.
\end{equation}
This definition of $\bs$ will be kept throughout the rest of this paper.

Let us define a sequence $(s_n)_{n\geq1-k}$ of words associated with $\bs$ as
follows:
\begin{align}
  &s_{1-k} = a_2, \quad s_{2-k} = a_{3}, \quad \ldots~, \quad
  s_{-1} = a_k, \quad s_{0} = a_1, \notag \\
  &s_n = s_{n-1}^{d_{n}}s_{n-2}^{d_{n-1}}\cdots
                           s_{0}^{d_{1}}a_{n+1}, \quad
   1\leq n \leq k-1, \label{eq:7.4} \\
  &s_{n} = s_{n-1}^{d_{n}}s_{n-2}^{d_{n-1}}\cdots
                           s_{n-k+1}^{d_{n-k+2}}s_{n-k}, \quad n \geq k.
   \notag
\end{align}
Clearly, $s_{n}$ is a prefix of $s_{n+1}$ for all $n \geq 0$ (and
hence $(|s_{n}|)_{n\geq0}$ is a strictly increasing sequence of
positive integers).


\begin{example} \label{ex:sturm&epis} 
It is well-known that the standard (or characteristic) Sturmian word $c_\alpha$ of 
irrational \emph{slope} $\alpha = [0; 1+d_1, d_2, d_3, \ldots]$, $d_1 \geq 1$, 
(see \cite{jBpS02stur} for a definition) is the standard episturmian word 
over $\cA = \{a,b\}$ with directive word $\Delta(c_\alpha) =
a^{d_1}b^{d_2}a^{d_3}b^{d_4}a^{d_5}\cdots$. 
We have $c_\alpha = \lim_{n\rightarrow\infty} s_n$, where $(s_n)_{n\geq-1}$ is 
the \emph{standard sequence} associated with $c_\alpha$, defined by
\[
  s_{-1} = b, \quad s_{0} = a, \quad s_n = s_{n-1}^{d_{n}}s_{n-2}, \quad 
  n\geq 1.
\]
This coincides with our definition \eqref{eq:7.4} above. 
Observe that, 
for all $n\geq 0$, $|s_n| = q_n$, where $q_n$ is the denominator of the 
$n$-th convergent to $[0; 1+d_1, d_2, d_3, \ldots]$.

For all $m \geq 1$, let $L_m := d_1 + d_2 + \cdots + d_m$. Then, writing 
$\Delta(c_\alpha) = x_1x_2x_3\cdots$ with each $x_i \in \cA$, we have 
$x_{n+1} \ne x_n$ if and only if $n$ is equal to some $L_m$. 
One easily deduces that $S(L_m) = L_{m+1} + 1$ and $P(L_{m+1} + 1) = L_m$, 
and it can also be shown that the $h_{L_m}$ satisfy the same recurrence 
relation as the $q_m$. Hence, $|h_{L_m}| = q_m$, and clearly we have
$h_{L_m} = s_m$ (see Proposition \ref{P:19.11.03} below).
\qed
\end{example}

Notice that $\bs$ has directive word resembling $\Delta(c_\alpha)$. 

\begin{notation}
Hereafter, let $L_n := d_1 + d_2 + \cdots + d_n$ for each
$n\geq1$.
\end{notation}

\begin{proposition} \label{P:19.11.03}
For any $n \geq 1$, $s_n = h_{L_n}$. Moreover, $\bs =
\underset{n\rightarrow\infty}{\lim}s_n$.
\end{proposition}
\begin{proof} 
The directive word of $\bs$ is given by
\[
  \Delta = a_1^{d_1}a_{2}^{d_2}\cdots a_k^{d_k}a_{1}^{d_{k+1}}
           a_{2}^{d_{k+2}}\cdots a_{k}^{d_{2k}}a_{1}^{d_{2k+1}}\cdots
         = x_1x_2x_3x_4\cdots, \quad x_i \in \cA_k.
\]

For $n \geq 1$, we have $x_{n+1} \ne x_n$ (and hence $h_{n} \ne
h_{n-1}$) if and only if $n$ is equal to some $L_m$. In
particular, for any $m \geq 1$,
\begin{equation} \label{eq:triangle_star}
  h_{L_m} = h_{L_{m+1} - r}, \quad 1 \leq r \leq d_{m+1}.
\end{equation}
Furthermore, it is clear that, for all $n \geq k$,
\begin{equation} \label{eq:triangle_star*}
  P(L_{n}+1) = L_{n-k+1}, 
\end{equation}
and $P(L_n+1)$ is undefined for $1\leq n\leq k-1$.

First we show that $s_n = h_{L_n}$ for $1 \leq n \leq k$. Observe
that, for $1 \leq n \leq k - 1$,
\begin{align*}
  h_{L_n} &= \Psi_{a_1}^{d_1}\Psi_{a_2}^{d_2}\cdots\Psi_{a_n}^{d_n}(a_{n+1}),\\
          &= \Psi_{a_1}^{d_1}\Psi_{a_2}^{d_2}\cdots\Psi_{a_{n-1}}^{d_{n-1}}
             (a_{n}^{d_n}a_{n+1}) \\
          &= h_{L_{n-1}}^{d_n}\Psi_{a_1}^{d_1}\Psi_{a_2}^{d_2}\cdots
             \Psi_{a_n}^{d_n}(a_{n+1}) \\
          &= h_{L_{n-1}}^{d_n}h_{L_{n-2}}^{d_{n-1}}\cdots h_{L_1}^{d_2}
             \Psi_{a_1}^{d_1}(a_{n+1}) \\
          &= h_{L_{n-1}}^{d_n}h_{L_{n-2}}^{d_{n-1}}\cdots h_{L_1}^{d_2}
             {a_1}^{d_1}a_{n+1} \\
          &= h_{L_{n-1}}^{d_n}h_{L_{n-2}}^{d_{n-1}}\cdots h_{L_1}^{d_2}
             h_0^{d_1}a_{n+1}.
\end{align*}
Similarly, since $h_{L_k} =
\Psi_{a_1}^{d_1}\Psi_{a_2}^{d_2}\cdots\Psi_{a_k}^{d_k}(a_{1})$,
one finds that
\[
  h_{L_k} = h_{L_{k-1}}^{d_k}h_{L_{k-2}}^{d_{k-1}}\cdots h_{L_1}^{d_2}a_{1}.
\]
Thus, we see that the $s_n$ satisfy the same recurrence relation
as the $h_{L_n}$ for $1 \leq n \leq k$. Therefore, since 
$h_0 = \mu_0(a_1) = a_1 = s_0$, we have
\begin{equation} \label{eq:15.11-1}
  s_n = h_{L_n} \quad \mbox{for all $n$, $1 \leq n \leq k$}.
\end{equation}

Now take $n \geq k+1$.  
Then, by \eqref{eq:h_n-1} and \eqref{eq:triangle_star*}, we have
\[
  h_{L_n} = h_{L_{n}-1}h_{L_{n}-2}\cdots h_{L_{n-k+1}}h_{L_{n-k+1}-1},
\]
and therefore it follows from \eqref{eq:triangle_star} that
\begin{equation} \label{eq:triangle_star**}
  h_{L_n} =  h_{L_{n-1}}^{d_{n}}h_{L_{n-2}}^{d_{n-2}}\cdots
             h_{L_{n-k+1}}^{d_{n-k+2}}h_{L_{n-k}}.
\end{equation}
Whence, since $s_n = h_{L_n}$ for $1 \leq n \leq k$,
\eqref{eq:triangle_star**} shows that the $s_n$ satisfy the same
recurrence relation as the $h_{L_n}$ for $n \geq k+1$. Thus, by virtue of 
this fact and \eqref{eq:15.11-1}, we have
\[
  s_n = h_{L_n} \qquad \mbox{for all  $n \geq 1$},
\]
as required.

The second assertion follows immediately from the first since $\bs
= \underset{m\rightarrow\infty}{\lim}h_m$ \cite{jJgP02epis}.
\end{proof}

Accordingly, the words $(s_n)_{n\geq1}$ can be viewed as `building blocks' of $\bs$.

\begin{example} The \emph{Tribonacci sequence} (or \emph{Rauzy word} \cite{gR82nomb}) is the standard 
episturmian word over $\{a,b,c\}$ directed by $(abc)^\omega$. Since 
all $d_i = 1$, we have $L_n = n$, and hence 
$h_n = s_n = s_{n-1}s_{n-2}s_{n-3}$  for all $n\geq1$. \qed
\end{example}

\subsection{Two special integer sequences}

Set $Q_n := |s_n|$ for all $n \geq 0$. Then the integer sequence
$(Q_n)_{n\geq0}$ is given by:
\begin{align*}
  Q_0 &= 1, \quad Q_n = d_nQ_{n-1} + d_{n-1}Q_{n-2} + \cdots + d_1Q_0 + 1,
  \quad 1\leq n \leq k-1, \\
  Q_n &= d_{n}Q_{n-1} + d_{n-1}Q_{n-2} + \cdots + d_{n+2-k}Q_{n+1-k} +
  Q_{n-k}, \quad n \geq k.
\end{align*}
Now, define the integer sequence $(P_n)_{n\geq0}$ by:
\begin{align*}
P_0 &= 0, \quad P_n = d_nP_{n-1} + d_{n-1}P_{n-2} + \cdots +
d_1P_0 + 1, \quad
          1 \leq n \leq k-1, \\
P_n &= d_nP_{n-1} + d_{n-1}P_{n-2} + \cdots +  d_{n+2-k}P_{n+1-k}
+
       P_{n-k}, \quad n \geq k.
\end{align*}

For $k = 2$, observe that $P_n/Q_n$ is the $n$-th convergent to
the continued fraction expansion $[0; 1+d_1, d_2, d_3, d_4,
\ldots]$.

\begin{proposition} \label{P:13.09.04(1)}
For all $n \geq 0$, $|s_n|_{a_1} = Q_n - P_n$.
\end{proposition}
\begin{proof} 
Proof is by induction on $n$. 
\end{proof}

\section{Generalized singular words} \label{S:g_singular}

Recall the standard Sturmian word $c_\alpha$ of slope 
$\alpha = [0;1+d_1,d_2,d_3,\ldots]$, $d_1 \geq 1$. Melan\c{c}on 
\cite{gM99lynd} (also see \cite{wCzW03some}) introduced the \emph{singular words} $(w_n)_{n\geq1}$ of 
$c_\alpha$ defined by
\[
  w_n = \begin{cases}
         as_nb^{-1} ~~\mbox{if $n$ is odd}, \\
         bs_na^{-1} ~~\mbox{if $n$ is even},
        \end{cases}
\]
with the convention $w_{-2} = \empt$, $w_{-1} = a$, $w_{0} = b$. It is easy to show that the set of factors of $c_\alpha$ of length $|s_n|$ is given by
\[
   \Omega_{|s_n|}(c_\alpha) = \cC(s_n) \cup \{w_n\}.
\]
(See \cite{gM99lynd, wCzW03some, aG04occu} for instance.) Also note that in this 2-letter case $s_n = u_{L_n}ab$ (resp.~$s_n = u_{L_n}ba$) if $n$ is odd (resp.~even).

Singular words are profoundly useful in studying properties of factors 
of $c_\alpha$ (e.g., \cite{wCzW03some, aG04occu, aG04conj, fLpS03conj, gM99lynd, zWzW94some}). 
It is for this very reason that we now generalize these 
words for the standard episturmian word $\bs$. Firstly, however, we 
prove some basic results concerning the words $s_n$ and $u_{L_{n}}$, as 
detailed in the next section.

\subsection{Useful results}

For each $n \geq 0$, set $D_{n} := u_{L_{n+1}}$. Observe that, for
any $m \geq 1$,
\begin{equation} \label{eq:starry1}
  |D_m| = (d_{m+1}-1)|s_m| + \sum_{j=0}^{m-1}d_{j+1}|s_j|.
\end{equation}
Indeed, using \eqref{eq:u_n} and \eqref{eq:triangle_star}, one finds that 
\begin{align}
D_m = u_{L_{m+1}}  &= h_{L_{m+1}-2}h_{L_{m+1}-3}\cdots h_{1}h_0 \notag \\
                       &= h_{L_{m}}^{d_{m+1}-1}h_{L_{m-1}}^{d_m}
                          h_{L_{m-2}}^{d_{m-1}}\cdots
                          h_{L_1}^{d_2}h_0^{d_1}  \notag \\
                       &= s_{m}^{d_{m+1}-1}s_{m-1}^{d_m}
                           s_{m-2}^{d_{m-1}}\cdots
                          s_1^{d_2}s_0^{d_1}. \label{eq:15.04.04}
\end{align}
Also note that $D_0 = a_1^{d_{1}-1}$ since $D_0 = u_{d_1} =
h_{d_1-2}h_{d_1-3}\cdots h_1h_0 = h_0^{d_1-1}$. For technical
reasons, we shall set $D_{-j} := a_{k+1-j}^{-1}$ and $|D_{-j}| =
-1$ for $1\leq j \leq k$.

\begin{proposition} \label{P:27.11.03(1)}
Let $1 \leq i \leq k$. For all $n \geq 1-k$, $a_i$ is the last letter of  $s_{n}$ if 
$n \equiv i-1 \pmod{k}$.
\end{proposition}
\begin{proof} 
Since we have $s_{1-k} = a_2$, $s_{2-k} = a_{3}$, 
$\ldots$~, $s_{-1} = a_k$, $s_{0} = a_1$,
the result follows immediately from the definition of the words
$s_n$ (see \eqref{eq:7.4}).
\end{proof}

\begin{proposition} \label{P:03.02.04(1)}
For all $n \geq 0$, $s_{n+1}D_{n-k+1} = s_nD_n$, and hence $|D_n|
- |D_{n-k+1}| = |s_{n+1}| - |s_n|$.
\end{proposition}
\begin{proof}
The claim holds for $0 \leq n \leq k-2$ since
$s_{n+1}D_{n-k+1} = s_n^{d_{n+1}}\cdots s_{0}^{d_1}a_{n+2}a_{n+2}^{-1} = 
s_nD_n$, and for $n\geq k-1$, 
$s_{n+1}D_{n-k+1} = s_{n}^{d_{n+1}}s_{n-1}^{d_{n}}\cdots s_{n-k+1}^{d_{n-k+2}}
                 \cdots s_{1}^{d_2}s_{0}^{d_1} = s_{n}D_{n}$.
\end{proof}

\begin{proposition} \label{P:02.09.04(1)}
For all $n \geq 1$, $|s_n| > |D_{n-1}|$. 
\end{proposition}
\begin{proof} Proof is by induction on $n$. The result is clearly true
for $n=1$ since $|s_1| = |a_1^{d_1}a_2| = |D_0a_1a_2| = |D_0| +
2$. Now assume the result holds for some $n \geq 2$. Then, using
Proposition \ref{P:03.02.04(1)},
\[
  |s_{n+1}| = |s_n| + |D_n| - |D_{n-k+1}|
           > |D_{n-1}| + |D_n| - |D_{n-k+1}|
            \geq |D_n|,
\]
since $|D_{n-k+1}| \leq |D_{n-1}|$.
\end{proof}

Recall that the words $D_n$ and $s_n$ are prefixes of $\bs$ for
all $n \in \NN$. Thus, according to Proposition
\ref{P:02.09.04(1)}, the palindromes $D_0$, $D_1$, $\ldots$~,
$D_{n-1}$ are prefixes of $s_n$. In fact, the maximal index $i$
such that $D_{i}$ is a proper prefix of $s_n$ is $i = n-1$, which
is evident from the following result.

\begin{proposition} \label{P:1.12.03(2)}
For all $n \geq 0$, $D_n = s_{n}^{d_{n+1}}D_{n-k}$.
\end{proposition}
\begin{proof}
Firstly, $D_0 = a_1^{d_1-1} = s_0^{d_{1}}a_{1}^{-1} =
s_0^{d_{1}}D_{-k}$ and, for $1 \leq n \leq k-1$, we have 
\begin{align*}
  D_{n} &= s_{n}^{d_{n+1}-1}s_{n-1}^{d_n}  \cdots s_0^{d_{1}} \\
             & = s_n^{d_{n+1}-1}s_na_{n+1}^{-1} \quad \mbox{(using \eqref{eq:7.4})} \\
             &= s_{n}^{d_{n+1}}a_{n+1}^{-1} = s_n^{d_{n+1}}D_{n-k}. 
  \end{align*}
Now take $n \geq k$. Then
\[
  D_n =  s_{n}^{d_{n+1}-1}s_{n-1}^{d_n}
  \cdots s_{n-k+1}^{d_{n-k+2}}s_{n-k}D_{n-k} = s_{n}^{d_{n+1}-1}s_nD_{n-k}
  = s_n^{d_{n+1}}D_{n-k}.
\]
\end{proof}

\begin{proposition} \label{P:05.02.04(1)-2}
For all $n \geq 0$, $s_n = D_{n-k}\rev{s}_nD_{n-k}^{-1}$.
\end{proposition}
\begin{proof}
Proof is by induction on $n$. For $n = 0$,
$D_{-k}\rev{s}_0D_{-k}^{-1} = a_{1}^{-1}a_1a_1 = a_1 = s_0$.
Assume the result holds for some $n \geq 1$. Then, using
Proposition \ref{P:03.02.04(1)},
\[
s_{n+1}
=s_nD_nD_{n-k+1}^{-1}
=D_{n-k}\rev{s}_nD_{n-k}^{-1}D_nD_{n-k+1}^{-1}.
\]
Therefore, invoking Proposition \ref{P:1.12.03(2)} and \eqref{eq:7.4}, for $1 \leq n
\leq k-2$, we have
\begin{align*}
  s_{n+1} &= D_{n-k}\rev{s}_n(\rev{s}_n)^{d_{n+1}}D_{n-k+1}^{-1} \\
                &= D_{n-k+1}a_{n+2}a_{n+1}^{-1}a_{n+1}(\rev{s}_0)^{d_1}\cdots
                   (\rev{s}_{n-1})^{d_{n}}(\rev{s}_n)^{d_{n+1}}
                   D_{n-k+1}^{-1} \\
                &= D_{n-k+1}a_{n+2}a_{n+2}^{-1}\rev{s}_{n+1}D_{n-k+1}^{-1} \\
                &= D_{n-k+1}\rev{s}_{n+1}D_{n-k+1}^{-1}.
\end{align*}
And, for $n\geq k-1$,
\begin{align*}
s_{n+1} &= D_{n-k}\rev{s}_n(\rev{s}_n)^{d_{n+1}}D_{n-k+1}^{-1} \\
        &= D_{n-k}[\rev{s}_{n-k}(\rev{s}_{n-k+1})^{d_{n-k+2}-1}
           \rev{s}_{n-k+1}(\rev{s}_{n-k+2})^{d_{n-k+3}}
           \cdots (\rev{s}_{n-1})^{d_{n}}](\rev{s}_{n})^{d_{n+1}}
            D_{n-k+1}^{-1} \\
        &= D_{n-k+1}\rev{s}_{n+1}D_{n-k+1}^{-1},
\end{align*}
as required.
\end{proof}

\begin{remark} \label{R:09.01.05}
This result shows, in particular, that $\rev{s}_n =
D_{n-k}^{-1}s_nD_{n-k}$, i.e., $\rev{s}_n$ is the $|D_{n-k}|$-th
conjugate of $s_n$ for each $n\geq k$. (For $0 \leq n \leq k-1$,
$\rev{s}_n$ is the $(|s_n|-1)$-st conjugate of $s_n$ since
$\rev{s}_n = a_{n+1}s_na_{n+1}^{-1}$.) The following two 
corollaries are direct results of the above proposition.
\end{remark}

\begin{corollary} \label{C:08.09.04(1)}
For any $n\geq 0$, the word $\rev{s}_nD_{n-k}^{-1}$ is a
palindrome. In particular, let $U_n= D_{n-k}$ and
$V_n=\rev{s}_nD_{n-k}^{-1}$. Then $s_n = U_nV_n$ is the unique
factorization of $s_n$ as a product of two palindromes.
\end{corollary}
\begin{proof}
From Proposition \ref{P:05.02.04(1)-2}, we have $s_n =
D_{n-k}\rev{s}_{n}D_{n-k}^{-1} = U_nV_n$, and whence $D_{n-k}^{-1}s_n
= \rev{s}_{n}D_{n-k}^{-1}$. It is therefore clear that
$\rev{s}_nD_{n-k}^{-1}$ is a palindrome. The uniqueness of the
factorization $s_n = U_nV_n$ is immediate from the primitivity of
$s_n$, which follows from Lemma \ref{L:h_n&u_n}(i), together with
Proposition \ref{P:19.11.03}. (Recall that since $s_n$ is primitive,
there are exactly $|s_n|$ different conjugates of $s_n$.)
\end{proof}

\begin{corollary} \label{C:05.02.04(1)}
For all $n \geq 0$, $s_n = D_{n}\rev{s}_nD_{n}^{-1}$.
\end{corollary}
\begin{proof} Propositions \ref{P:1.12.03(2)} and \ref{P:05.02.04(1)-2}. 
\end{proof}

\begin{notation}
Now, for each $n \in \NN$, we define the words $G_{n,r}$ by
\[
  s_n = D_{n-r}G_{n,r}, \quad 1\leq r\leq k-1.
\]
\end{notation}

\begin{example}
In the case of Sturmian words $c_\alpha$, $r = 1$
and $s_n = D_{n-1}G_{n,1} = u_{L_n}G_{n,1}$ for all $n\geq1$, 
where $G_{n,1} = ab$ or $ba$, according to $n$ odd or even, respectively. \qed
\end{example}

\begin{example} Recall that when all $d_i=1$, $\bs$ is the Tribonacci sequence over $\{a_1,a_2,a_3\} \equiv \{a,b,c\}$. For $n = 4$,  we have $s_n =s_4 = abacabaabacab$, $D_2 = aba$, $D_3 = abacaba$, and hence 
\[
  G_{4,1} = abacab \quad \mbox{and} \quad G_{4,2} = cabaabacab. 
\]  
\qed
\end{example}

\begin{note} Since $D_{n-r} = a_{k+1+n-r}^{-1}$ for $0 \leq n
< r$, we also set 
\begin{equation} \label{eq:22.09-1}
G_{n,r} = a_{k+1+n-r}s_n \quad \mbox{for  $0\leq n < r$}.
\end{equation}
\end{note}

\begin{proposition} \label{P:12.05.04(1)} 
For all $n \geq 1$, $s_{n}s_{n-1}G_{n-1,k-1}^{-1} =
s_{n-1}s_{n}G_{n,1}^{-1}$.
\end{proposition}
\begin{proof} 
It is easily checked that the result holds for $1\leq n \leq k-1$,
since
\[
s_ns_{n-1}G_{n-1,k-1}^{-1} = s_nD_{n-k} = s_na_{n+1}^{-1},
\]
and
\[
  s_{n-1}s_nG_{n,1}^{-1} = s_{n-1}D_{n-1} = s_{n-1}^{d_{n}}\cdots s_0^{d_1}
                         = s_{n}a_{n+1}^{-1}.
\]
Now take $n\geq k$. Then, using \eqref{eq:15.04.04}, we have
\begin{align*}
  s_{n}s_{n-1}G_{n-1,k-1}^{-1} &= s_{n}D_{n-k} \\
  &= (s_{n-1}^{d_{n}}s_{n-2}^{d_{n-1}}\cdots
     s_{n-k+1}^{d_{n-k+2}}s_{n-k})s_{n-k}^{d_{n-k+1}-1}s_{n-k-1}^{d_{n-k}}
     \cdots s_{1}^{d_2}s_{0}^{d_1} \\
  &= s_{n-1}(s_{n-1}^{d_{n}-1}s_{n-2}^{d_{n-1}}\cdots
     s_{n-k+1}^{d_{n-k+2}}s_{n-k}^{d_{n-k+1}}s_{n-k-1}^{d_{n-k}}
     \cdots s_{1}^{d_2}s_{0}^{d_1}) \\
  &= s_{n-1}D_{n-1} \\
  &= s_{n-1}s_{n}G_{n,1}^{-1}.
\end{align*}
\end{proof}

\begin{remark} \label{R:13.05.04}
Recall Example \ref{ex:sturm&epis}. For $c_\alpha$ with $\alpha =
[0;1+d_1,d_2,d_3\ldots]$, it is well-known that, for all $n \geq
2$, $s_ns_{n-1}(xy)^{-1} = s_{n-1}s_n(yx)^{-1}$, where $x,y\in
\{a,b\}$, $x\ne y$, and $xy \suff s_{n-1}$. This is known as the
\emph{Near-Commutative Property} of the words $s_n$ and $s_{n-1}$.
Because $s_ns_{n-1}(xy)^{-1} = s_{n}D_{n-2}$ and
$s_{n-1}s_{n}(yx)^{-1} = s_{n-1}D_{n-1}$, Proposition
\ref{P:12.05.04(1)} is merely an extension of this property to 
standard episturmian words $\bs$. It is also worthwhile noting that
Proposition \ref{P:12.05.04(1)} shows that $s_n$ is a prefix of
$s_{n-1}s_n$. 
\end{remark}

Hereafter, we set $d_{-j} = 0$ for $j \geq 0$.

Proposition \ref{P:03.02.04(1)} implies that $|s_{n+1}| - |D_n| =
|s_n| - |D_{n-k+1}|$, and hence $|G_{n+1,1}| = |G_{n,k-1}|$. In
fact, we have the following:

\begin{proposition} \label{P:G-wiz}
For all $n \geq 1$, $G_{n,1} = \rev{G}_{n-1,k-1}$.
\end{proposition}
\begin{proof} 
One can write
\begin{align*}
G_{n,1} = D_{n-1}^{-1}s_{n} &=
D_{n-1}^{-1}s_{n-1}^{d_{n}}s_{n-2}^{d_{n-1}}\cdots
                         s_{n-k+1}^{d_{n-k+2}}s_{n-k} \\
          &= D_{n-1}^{-1}s_{n-1}s_{n-1}^{d_{n}-1}s_{n-2}^{d_{n-1}}\cdots
                         s_{n-k+1}^{d_{n-k+2}}s_{n-k} \\
          &= D_{n-1}^{-1}s_{n-1}D_{n-1}D_{n-k}^{-1}.
\end{align*}
Whence, it follows from Corollary \ref{C:05.02.04(1)}
that $G_{n,1} = \rev{s}_{n-1}D_{n-k}^{-1} = \rev{G}_{n-1,k-1}$ 
since $\rev{s}_{n-1} = \rev{G}_{n-1,k-1}D_{n-k}$.
\end{proof}

\newpage
\begin{proposition} \label{P:19.04.04(2)}
Let $1 \leq i \leq k$ and $1\leq r \leq k-1$. For all $n\geq 0$,
\begin{itemize}
\item[\emph{(i)}] ~$a_i$ is the first letter of $G_{n,r}$ if $n \equiv i + r - 1 \pmod{k};$
\item[\emph{(ii)}] ~$a_i$ is the last letter of $G_{n,r}$ if $n \equiv i - 1 \pmod{k}$.
\end{itemize}
\end{proposition}
\begin{proof} 
(i) The assertion is trivially true for $0 \leq n < r$ since,  by 
\eqref{eq:22.09-1},  we have $G_{n,r} = a_{k+1+n-r}s_n$.  
Now take $n \geq r$. By definition,
\[
  G_{n,r} = D_{n-r}^{-1}s_{n} = D_{n-r}^{-1}s_{n-r+1}s_{n-r+1}^{-1}s_n 
\]
where $s_{n-r+1}$ is a prefix of $s_n$. Hence, one can write
\begin{equation} \label{eq:20.09.04(2)}
G_{n,r} = G_{n-r+1,1}s_{n-r+1}^{-1}s_{n}
   = \rev{G}_{n-r,k-1}s_{n-r+1}^{-1}s_{n} 
\end{equation}
by applying Proposition \ref{P:G-wiz}.

Now, one easily deduces from Proposition \ref{P:27.11.03(1)} that
$a_m \pref \rev{s}_{n-r}$ if $n \equiv m+r-1 \pmod{k}$, and thus
$a_m \pref \rev{G}_{n-r,k-1} \pref \rev{s}_{n-r}$ if $n \equiv
m+r-1 \pmod{k}$.

(ii) For $0 \leq n < r$, $G_{n,r} = a_{k+1+n-r}s_n$ and, for each
$n \geq r$, we have $G_{n,r} \suff s_{n}$. Hence, $a_m \suff
G_{n,r}$ if $n \equiv m-1 \pmod{k}$, by Proposition
\ref{P:27.11.03(1)}.
\end{proof}

\subsection{Singular $n$-words of the $r$-th kind}

By definition of the words $(s_n)_{n\geq1-k}$ (see \eqref{eq:7.4})
and the fact that $\bs = \lim_{n\rightarrow\infty}s_n$, one
deduces that, for any $n\geq0$, $\bs$ can be written as a
concatenation of blocks of the form $s_n$, $s_{n-1}$, $\ldots$~,
$s_{n-k+1}$, i.e.,
\begin{align}
  \bs = &~[((s_{n}^{d_{n+1}}s_{n-1}^{d_n}\cdots s_{n-k+2}^{d_{n-k+3}}
        s_{n-k+1})^{d_{n+2}}s_n^{d_{n+1}}\cdots s_{n-k+3}^{d_{n-k+4}}
        s_{n-k+2})^{d_{n+3}} \notag \\ &~(s_{n}^{d_{n+1}}s_{n-1}^{d_n}\cdots
        s_{n-k+2}^{d_{n-k+3}}s_{n-k+1})^{d_{n+2}}s_{n}^{d_{n+1}}\cdots
        s_{n-k+4}^{d_{n-k+5}}s_{n-k+3}]^{d_{n+4}}\cdots. \label{eq:n-part-k}
\end{align}
We shall call this unique decomposition the
\emph{$n$-partition} of $\bs$. This will be a useful tool in our
subsequent analysis of powers of words occurring in $\bs$ (Section
\ref{S:powers_in_s}, to follow).

\begin{note} Uniqueness of the factorization \eqref{eq:n-part-k} is proved
inductively. The initial case $n=0$ is trivial. For $n\geq 1$, the factorization of $s_n$ in terms of the 
$s_{n-i}$ given by \eqref{eq:7.4} is unique because the $s_{n-i}$ end with different letters (by Proposition \ref{P:27.11.03(1)}). So it is clear that every $(n+1)$-partition of $\bs$ gives rise to an
$n$-partition, in which the positions of $s_{n-k+1}$ blocks
uniquely determine the positions of $s_{n+1}$ blocks in the
original $(n+1)$-partition (since $s_{n+1} =
s_n^{d_{n+1}}s_{n-1}^{d_{n}}\cdots s_{n-k+2}^{d_{n-k+3}}
s_{n-k+1}$). Accordingly, uniqueness of the $n$-partition implies uniqueness of the $(n+1)$-partition.
\end{note}

\begin{remark} \label{R:09.09.04}
Since each factor of $\bs$ has exactly $k$ different return words (see Section \ref{SS:return}), 
two consecutive $s_{n+1-i}$ blocks ($1\leq i \leq k$) of the $n$-partition 
are separated by a word $V$, of which there are $k$ different possibilities. 
From now on, it is advisable to keep this observation in mind.
\end{remark}

\begin{lemma} \label{L:06.09.04(1)} Let $1 \leq r \leq k-1$.
For any $n\in \NN^+$, a factor $u$ of length $|s_n|$ of $\bs$ is a
factor of at least one of the following words:
\begin{itemize}
\item $C_j(s_n),$ $0\leq j\leq |s_n| - 1;$
\item $s_{n-r}^{d_{n-r+1}-1}\cdots s_{n-k+1}^{d_{n-k+2}}s_{n-k}
       s_{n-1}^{d_{n}}\cdots s_{n-r+1}^{d_{n-r+2}}s_{n-r}s_n \quad
      \mbox{if ~$n\geq r;$}$
\item $a_{n+1}s_na_{n+1}^{-1}a_{n-r+k+1}s_n \quad \mbox{if ~$n < r$}$.
\end{itemize}
\end{lemma}

\begin{note} 
The word $s_{n-r}^{d_{n-r+1}-1}\cdots
s_{n-k+1}^{d_{n-k+2}}s_{n-k}
       s_{n-1}^{d_{n}}\cdots s_{n-r+1}^{d_{n-r+2}}s_{n-r}$ ($1\leq r\leq k-1$) 
has length $|s_n|$.
\end{note}

\begin{proof}[Proof of Lemma $\ref{L:06.09.04(1)}$] 
In the $n$-partition of $\bs$, one observes that two consecutive
$s_n$ blocks make the following $k$ different appearances:
\[
  s_ns_n \quad \mbox{and} \quad
  \underbrace{s_ns_{n-1}^{d_{n}}\cdots s_{n-r+1}^{d_{n-r+2}}
  s_{n-r}s_n}_{(*)}, \quad 1 \leq r \leq k-1.
\]
Evidently, any factor of length $|s_n|$ of $\bs$ is a factor of one
of the above $k$ different words. 

Now, factors of length $|s_n|$ of $s_ns_n$ are simply conjugates
of $s_n$. Furthermore, for $n\geq r$, the first
$|s_{n-1}^{d_{n}}\cdots s_{n-r+1}^{d_{n-r+2}}s_{n-r}|$ factors of
length $|s_n|$ of $(*)$ are again just conjugates of $s_n$. The
remaining factors of length $|s_n|$ of $(*)$ are factors of
\[
  s_{n-r}^{d_{n-r+1}-1}\cdots s_{n-k+1}^{d_{n-k+2}}s_{n-k}
       s_{n-1}^{d_{n}}\cdots s_{n-r+1}^{d_{n-r+2}}s_{n-r}s_n.
\]
For $n < r$, one can write $(*)$ as
$s_ns_{n-1}^{d_n}\cdots s_0^{d_1}a_{n-r+k+1}s_n = s_ns_na_{n+1}^{-1}
a_{n-r+k+1}s_n$,
of which the first $|s_n|-1$ factors of length $|s_n|$ are
conjugates of $s_n$, and the other factors of length $|s_n|$ are
factors of $a_{n+1}s_na_{n+1}^{-1}a_{n-r+k+1}s_n$.
\end{proof}

\begin{lemma} \label{L:01.12.03(1)}
For any $n \geq 1$, $\sum_{j=1}^{k-1}|D_{n-j}| = |s_n| - k$.
\end{lemma}
\begin{proof}  
Induction on $n$ and Proposition \ref{P:03.02.04(1)}.
\end{proof}

\begin{lemma} \label{L:06.09.04(2)} Let $1\leq r \leq k-1$.
For any $n \geq r$, we have
\[
  s_{n-r}^{d_{n-r+1}-1}\cdots s_{n-k+1}^{d_{n-k+2}}s_{n-k}
       s_{n-1}^{d_{n}}\cdots s_{n-r+1}^{d_{n-r+2}}s_{n-r} =
  D_{n-r}\rev{G}_{n,r},
\]
and for $1\leq n < r$, $a_{n+1}s_na_{n+1}^{-1}a_{n-r+k+1} =
\rev{G}_{n,r}$.
\end{lemma}
\begin{proof} 
For $1 \leq n < r$, one can write
$\rev{G}_{n,r} = \rev{s}_na_{n-r+k+1} = a_{n+1}s_na_{n+1}^{-1}a_{n-r+k+1}$,
by Remark \ref{R:09.01.05}. 
Now take $n\geq r$. Then, using Corollary \ref{C:05.02.04(1)} and
Proposition \ref{P:1.12.03(2)},
\begin{align*}
D_{n-r}\rev{G}_{n,r} &= D_{n-r}\rev{s}_nD_{n-r}^{-1} \\
                     &= D_{n-r}D_{n}^{-1}s_nD_nD_{n-r}^{-1} \\
                     &= D_{n-r}D_n^{-1}s_n^{d_{n+1}}s_{n-1}^{d_{n}}\cdots
                        s_{n-r+1}^{d_{n-r+2}}s_{n-r} \\
                     &= D_{n-r}D_{n-k}^{-1}s_{n-1}^{d_{n}}\cdots
                        s_{n-r+1}^{d_{n-r+2}}s_{n-r} \\
                     &= s_{n-r}^{d_{n-r+1}-1}\cdots s_{n-k+1}^{d_{n-k+2}}
                        s_{n-k}s_{n-1}^{d_{n}}\cdots
                        s_{n-r+1}^{d_{n-r+2}}s_{n-r}.
\end{align*}
\end{proof}

Whence, it is now plain to see that each word 
$\rev{G}_{n,r}s_n = \rev{G}_{n,r}D_{n-r}G_{n,r}$ is a factor of $\bs$. We will
now partition the set of factors of length $|s_n|$ of $\bs$ into
$k$ disjoint classes.

\begin{theorem} \label{T:06.09.04(3)} Let $1\leq r\leq k-1$.
For any $n \in \NN^+$, the set of factors of length $|s_n|$ of $\bs$ can be partitioned
into the following $k$ disjoint classes:
\begin{itemize}
\item $\Omega_n^0 := \cC(s_n) = \{C_j(s_n) ~: ~ 0\leq j\leq |s_n|-1\};$
\item $\Omega_n^r := \{w \in \cA_k^* ~: ~ |w| = |s_n| ~\mbox{and} ~w \prec
      x^{-1}\rev{G}_{n,r}D_{n-r}G_{n,r}x^{-1}\}$,
where $x$ is the last letter of $G_{n,r}$.
\end{itemize}
That is, $\Omega_{|s_n|}(\bs) = \Omega_n^0\overset{\centerdot}{\cup}
\Omega_n^1\overset{\centerdot}{\cup}\cdots\overset{\centerdot}{\cup}\Omega_n^{k-1}$.
\end{theorem}
\begin{proof} 
First observe that Lemma \ref{L:h_n&u_n}(i), coupled with 
Proposition \ref{P:19.11.03}, implies that each $s_n$ is
primitive, and hence $|\Omega_n^0| = |s_n|$. Also note that
$\rev{\Omega}_n^0 := \{\rev{w} ~: ~ w \in \Omega_n^0\} =
\Omega_n^0$, i.e., $\Omega_n^0$ is closed under reversal, which is
deduced from Corollary \ref{C:08.09.04(1)}.

We shall use Lemma \ref{L:06.09.04(1)} to partition
$\Omega_{|s_n|}(\bs)$ into $k$ disjoint classes; the first being
$\Omega_n^0 = \cC(s_n)$. Now consider the factors of length $|s_n|$ of the words
\begin{equation}\label{eq:06.09.04-1}
  s_{n-r}^{d_{n-r+1}-1}\cdots s_{n-k+1}^{d_{n-k+2}}s_{n-k}
       s_{n-1}^{d_{n}}\cdots s_{n-r+1}^{d_{n-r+2}}s_{n-r}s_n  \quad (n\geq r). 
\end{equation}
Since \eqref{eq:06.09.04-1} can be written as
$D_{n-r}\rev{G}_{n,r}D_{n-r}G_{n,r}$ (by Lemma
\ref{L:06.09.04(2)}), the first $|D_{n-r}| + 1$ factors of length
$|s_n| = |D_{n-r}G_{n,r}|$ are conjugates of $\rev{s}_n$ (and
hence of $s_n$) and the last factor is just $s_n$. Hence, all
other factors of length $|s_n|$ of \eqref{eq:06.09.04-1} are
factors of $x^{-1}\rev{G}_{n,r}D_{n-r}G_{n,r}x^{-1}$, where $x$ is
the last letter of $G_{n,r}$. Moreover,  
$D_{n-r}$ appears exactly once (and at a different position) in
each word in $$\Omega_n^r := \{w \in \cA_k^* ~: ~ |w| = |s_n|
~\mbox{and} ~w \prec x^{-1}\rev{G}_{n,r}D_{n-r}G_{n,r}x^{-1}\};$$
whence $|\Omega_n^r| = |G_{n,r}| -1$. 
Since the letter just before $D_{n-r}$ (equivalently, the last letter
of $\rev{G}_{n,r}$) in the word
$x^{-1}\rev{G}_{n,r}D_{n-r}G_{n,r}x^{-1}$  is different for each $r \in [1,k-1]$, it is
evident that $\Omega_n^0$, $\Omega_n^1$, $\ldots$~,
$\Omega_n^{k-1}$ are pairwise disjoint.

Now, for $1\leq n < r$, other than words in the sets $\Omega_n^0$,
$\Omega_n^1$, $\ldots$~, $\Omega_n^n$, the remaining factors of
length $|s_n|$ of $\bs$ are factors of
\begin{equation} \label{eq:25.09.04}
  a_{n+1}s_na_{n+1}^{-1}a_{n-r+k+1}s_n = \rev{s}_na_{n-r+k+1}s_n
\end{equation}
(see Lemma \ref{L:06.09.04(1)}). The first factor of length
$|s_n|$ of the word \eqref{eq:25.09.04} is $\rev{s}_n$ (i.e., the
$(|s_n|-1)$-st conjugate of $s_n$) and the last is just $s_n$. All
other factors of length $|s_n|$ of  \eqref{eq:25.09.04} are
factors of
\[
  a_{n+1}^{-1}\rev{s}_na_{n-r+k+1}s_na_{n+1}^{-1}
  = a_{n+1}^{-1}\rev{G}_{n,r}D_{n-r}G_{n,r}a_{n+1}^{-1}.
\]
Defining $\Omega_n^r := \{ w\in \cA_k^* ~: ~ |w| = |s_n|
~\mbox{and} ~w \prec
a_{n+1}^{-1}\rev{G}_{n,r}D_{n-r}G_{n,r}a_{n+1}^{-1}\}$, one can
check
that $|\Omega_n^r| = |G_{n,r}| - 1$  
and $\Omega_n^0$, $\Omega_n^1$, $\ldots$~, $\Omega_n^{k-1}$ are
pairwise disjoint.

It remains to show $\bigcup_{j=0}^{k-1} \Omega_n^j =
\Omega_{|s_n|}(\bs)$ for all $n\geq 1$. Indeed,
$|\Omega_{|s_n|}(\bs)| = (k-1)|s_n| + 1$ (from the complexity
function for $k$-strict standard episturmian words), and we have
\begin{align*}
  \sum_{j=0}^{k-1}|\Omega_n^j| = |s_n| + \sum_{j=1}^{k-1}
      (|G_{n,j}| -1) 
  &= |s_n| + \sum_{j=1}^{k-1}(|s_n| - |D_{n-j}| - 1) \\
  &= k|s_n| - k + 1 - \sum_{j=1}^{k-1}|D_{n-j}| \\
  &= k|s_n| - k + 1 - (|s_n| - k) \quad \quad
     \mbox{(by Lemma \ref{L:01.12.03(1)})} \\
  &= (k-1)|s_n| + 1. 
\end{align*}
\end{proof}

Let us remark that the sets $\Omega_n^r$ are closed under reversal since $x^{-1}\rev{G}_{n,r}D_{n-r}G_{n,r}x^{-1}$ is a palindrome; that is $\rev{\Omega}_n^r := \{\rev{w} ~: ~w \in
\Omega_n^r\} = \Omega_n^r$. We 
shall call the factors of $\bs$ in $\Omega_n^r$ the \emph{singular
$n$-words of the $r$-th kind}. Such words will play a key role in
our study of powers of words occurring in $\bs$.

Evidently, for Sturmian words $c_\alpha$, $\Omega_n^1 = \{w_n\}$ and we have 
$\Omega_{|s_n|}(c_\alpha) = \cC(s_n)\cup\{w_n\}$, as before. 

\section{Index} \label{S:index}

A word of the form $w = (uv)^nu$ is written as $w = z^r$, where $z =
uv$ and $r := n + |u|/|z|$. The rational number $r$ is called the \emph{exponent} of $z$, and 
$w$ is said to be a \emph{fractional power}.

Now suppose $\bx$ is an infinite word. For any $w \prec \bx$, the
\emph{index} of $w$ in $\bx$ is given by the number
\[
  \mbox{ind}(w) = \sup\{r \in \QQ ~: ~w^r \prec \bx\},
\]
if such a number exists; otherwise, $w$ is said to have infinite
index in $\bx$. Furthermore, the greatest number $r$ such that
$w^r$ is a prefix of $\bx$ is called the \emph{prefix index} of
$w$ in $\bx$. Obviously, the prefix index is zero if the first
letter of $w$ differs from that of $\bx$, and it is infinite if
and only if $\bx$ is purely periodic.

For all $n \geq 0$, define the words
\[
  t_n := D_{n-k+1}G_{n+1,k-1} \quad \mbox{and} \quad
  r_n := s_{n-1}D_{n-1} = s_{n-1}^{d_{n}}s_{n-2}^{d_{n-1}}\cdots s_{1}^{d_{2}}
                          s_{0}^{d_{1}}.
\]

\begin{note} By convention, $r_0 = a_ka_k^{-1} = \empt$, and
$t_n = a_{n+2}^{-1}a_{n+3}s_{n+1}$ for $0 \leq n \leq k-2$.
\end{note}

The next two results extend those of Berstel \cite{jB99onth}.

\begin{lemma} \label{L:index1}
For all $n \geq 1$, the word $r_{n+1}$ is the greatest fractional
power of $s_{n}$ that is a prefix of $\bs$, and the prefix index
of $s_n$ in $\bs$ is $1 + d_{n+1} + |D_{n-k}|/|s_{n}|$.
\end{lemma}
\begin{proof} 
First we take $n\geq k$. Observe that the longest common prefix
shared by the words $s_n$ and $t_n$ is
\[
 D_{n-k+1} = s_{n-k+1}^{d_{n-k+2}-1}r_{n-k+1},
\]
since
\begin{equation} \label{eq:triangle-1}
s_{n} = D_{n-k+1}G_{n,k-1} \quad \mbox{and} \quad t_{n} =
D_{n-k+1}G_{n+1,k-1} = D_{n-k+1}\rev{G}_{n+2,1} 
\end{equation}
where $G_{n,k-1}$ and $G_{n+1,k-1}$  
do not share a common first letter, by Proposition
\ref{P:19.04.04(2)}. Clearly, $s_{n+1}s_n \pref \bs$, and we have
\begin{align}
  s_{n+1}s_n &= s_{n}^{d_{n+1}}s_{n-1}^{d_n}\cdots s_{n-k+2}^{d_{n-k+3}}
               s_{n-k+1}s_n \notag \\
             &= s_{n}^{d_{n+1}+1}(s_{n-k+1}^{d_{n-k+2}-1}s_{n-k})^{-1}
                D_{n-k+1}G_{n,k-1} \notag \\
              &= s_{n}^{d_{n+1}+1}(s_{n-k}^{d_{n-k+1}-1}s_{n-k-1}\cdots s_1^{d_2}s_0^{d_1})G_{n,k-1} 
                 \quad \mbox{(by \eqref{eq:15.04.04})}   \notag \\
             &= s_{n}^{d_{n+1}+1}D_{n-k}G_{n,k-1} \notag\\
             &=s_{n}^{d_{n+1}+1}t_{n-1}. \label{eq:triangle-2}
\end{align}
Hence, $s_{n}^{d_{n+1}+1}$ is a prefix of $\bs$. Also observe that
the longest common prefix of $t_{n-1}$ and $s_n$ is $D_{n-k}$ since 
\[
t_{n-1} = D_{n-k}G_{n,k-1} \quad \mbox{and} \quad  s_n = D_{n-k}G_{n,k}
\]
where $G_{n,k-1}$ and $G_{n,k}$ have different first letters, by Proposition \ref{P:19.04.04(2)}. Further, from \eqref{eq:triangle-2} and Proposition
\ref{P:1.12.03(2)}, we have
\[
 s_{n+1}s_{n} = s_{n}^{d_{n+1}+1}t_{n-1} = s_{n}^{d_{n+1}+1}D_{n-k}G_{n,k-1} =
              s_nD_nG_{n,k-1} =  r_{n+1}G_{n,k-1}.
\]
Thus, the greatest fractional power of $s_n$  that is a prefix of $\bs$
is $r_{n+1}$ with  
\[
  |r_{n+1}| = |s_{n}D_{n}| = |s_{n}^{d_{n+1}+1}D_{n-k}|
           = (d_{n+1}+1)|s_{n}| + |D_{n-k}|;
\]
whence the prefix index of $s_{n}$ in $\bs$ is $1 + d_{n+1} +
|D_{n-k}|/|s_n|$.

Similarly, for $1 \leq n \leq k-1$, we have 
\begin{align*}
  s_{n+1}s_n &= s_n^{d_{n+1}}s_{n-1}^{d_n}\cdots s_0^{d_1}a_{n+2}s_n \\
             &= s_n^{d_{n+1}+1}a_{n+1}^{-1}a_{n+2}s_n \\
             &= s_n^{d_{n+1}+1}D_{n-k}G_{n,k-1} \\
             &= r_{n+1}G_{n,k-1}. 
\end{align*}
Therefore, the greatest fractional power of $s_n$ (= $D_{n-k}G_{n-1,k-1}$)
 that is a prefix of $s_{n+1}s_n \pref \bs$
is $r_{n+1}$, where $|r_{n+1}| = (d_{n+1}+1)|s_n| + |D_{n-k}| =
(d_{n+1}+1)|s_n| - 1$. That is, the prefix index of $s_n$ in
$\bs$ is $1 + d_{n+1} - 1/|s_n|$ for $1\leq n \leq k-1$.
\end{proof}

\begin{lemma} \label{L:29.06.04(2)}
For all $n\geq 1$, the index of $s_n$ as a factor of $\bs$ is at
least $2 + d_{n+1} + |D_{n-k}|/|s_n|$, and hence $\bs$ contains
cubes.
\end{lemma}

We will show later that the index of $s_n$ is exactly $2 + d_{n+1}
+ |D_{n-k}|/|s_n|$.

\begin{proof} 
Setting $e = 1 + d_{n+1} + |D_{n-k}|/|s_n|$, we will show that
$s_{n+k+2}$ contains a power of $s_n$ of exponent $1+e$.
Certainly, using Proposition \ref{P:12.05.04(1)}, one can write
\begin{align*}
s_{n+k+2} &= s_{n+k+1}^{d_{n+k+2}-1}s_{n+k+1}s_{n+k}D_{n+k}D_{n+2}^{-1} \\
          &= s_{n+k+1}^{d_{n+k+2}-1}s_{n+k}s_{n+k+1}G_{n+k+1,1}^{-1}
            G_{n+k,k-1}D_{n+k}D_{n+2}^{-1} \\
          &= s_{n+k+1}^{d_{n+k+2}-1}s_{n+k}D_{n+k}G_{n+k,k-1}D_{n+k}
            D_{n+2}^{-1}.
\end{align*}
The suffix $s_{n+k}D_{n+k}G_{n+k,k-1}D_{n+k}D_{n+2}^{-1}$ contains
the exponent $1+e$ of $s_n$. More precisely, $s_{n+k}$ ends with
$s_{n}$, and $D_{n+k}G_{n+k,k-1}$ shares a prefix of length
$|D_{n+k}|$ with $s_{n+k+1}$. Thus, since $r_{n+1}$ is a prefix of
$\bs$ of length
\[
  |r_{n+1}| = |s_{n}| + |D_n| < |D_{n+k}|,
\]
we have $s_nr_{n+1} \prec s_{n+k}D_{n+k} \prec s_{n+k+2}$.
\end{proof}

\section{Powers occurring in $\bs$} \label{S:powers_in_s}

For each $m$, $l \in \NN$ with $l \geq 2$, let us define the
following set of words:
\[
  \cP(m;l) := \{w\in \cA_k^*~: ~|w| = m, ~w^l \prec \bs\},
\]
where $\bs$ is the $k$-strict standard episturmian word over
$\cA_k = \{a_1,a_2,\ldots, a_k \}$ with directive word $\Delta$
given by \eqref{eq:dir_seq}. Also, let $p(m;l) := |\cP(m;l)|$.

The next theorem is a generalization of Theorem 1 in
\cite{dDdL03powe}. It gives all the lengths $m$ such that there is
a non-trivial power of a word of length $m$ in $\bs$. Firstly, let us 
define the following $k$ sets of lengths for fixed $n \in
\NN^+$:
\begin{align*}
\cD_1(n) &:= \{r|s_n| ~: ~1\leq r \leq d_{n+1}\}, \\
\cD_i(n) &:= \{|s_n^rs_{n-1}^{d_{n}}\cdots
s_{n+2-i}^{d_{n+3-i}}s_{n+1-i}| ~:
            ~1\leq r \leq d_{n+1}\}, \quad 2 \leq i \leq k-1, \\
\cD_k(n) &:= \{|s_n^rs_{n-1}^{d_{n}}\cdots
s_{n+2-k}^{d_{n+3-k}}s_{n+1-k}|
            ~: ~1\leq r\leq d_{n+1} - 1\}.
\end{align*}

\begin{theorem} \label{T:22.11.03(5)}
Let $m$, $n \in \NN^+$  be such that $|s_n| \leq m < |s_{n+1}|$ and suppose $m \not\in
\bigcup_{i=1}^{k} \cD_i(n)$. Then 
$p(m;l) = 0$ for all $l\geq2$.
\end{theorem}

\begin{remark} Put simply, the above theorem states that if a word $w$ has a non-trivial integer power in $\bs$, then $|w| \in \bigcup_{i=1}^{k}\cD_i(n)$ for some $n$.   
For instance, if $k=3$, we have 
\[
  \bigcup_{i=1}^{3} \cD_i(n) = \{|s_n^r|, |s_{n}^rs_{n-1}| : 1 \leq r \leq d_{n+1}\}
  \cup \{|s_n^rs_{n-1}^{d_n}s_{n-2}| : 1 \leq r \leq d_{n+1}-1\}.
\]
In the particular case of the Tribonacci sequence, Theorem \ref{T:22.11.03(5)} implies that if $w^l$ is a factor, then 
$|w| \in \{|s_n|, |s_{n}| + |s_{n-1}|\}$ for some $n$, where the lengths $(|s_i|)_{i\geq0}$ are the  \emph{Tribonacci numbers}:  $T_0 = 1$, $T_1=2$, $T_2=4$, $T_i = T_{i-1} + T_{i-2} + T_{i-3}$, $i\geq 3$. 
\end{remark}

The proof of Theorem \ref{T:22.11.03(5)} requires several lemmas.
Let us first observe that in the $n$-partition of $\bs$ (see
\eqref{eq:n-part-k}) to the left of each $s_n$ block, there is an
$s_{n+1-j}$ block for some $j \in [1,k]$. Also note that each
$s_{n+1-j}$ is a prefix of $s_n$. Furthermore, to the left of each $s_{n+1-i}$ block is
another $s_{n+1-i}$ block or an $s_{n+2-i}$ block, for each $i \in
[2,k]$. 

\begin{lemma} \label{L:26.07.04(1)}
Let $n \in \NN^+$. 
Consider a word $w \prec \bs$ of the form $w = us_nv$ for some
words $u$, $v \in \cA_k^*$, $u \ne \empt$.
\begin{itemize}
\item[\emph{(i)}] If $w = u_1u_2$, where $u_1 \suff s_{n+1-i}$ for
some $i \in [1,k]$ and $u_2 \pref s_n$, then $u_1 = u$.
\item[\emph{(ii)}] If $w = u_1s_{n+1-i}u_2$ for some $i \in [2,k]$,
where $u_1 \suff s_{n+2-i}$ and  $u_2 \pref s_n$, then $u_1 = u$
or $u_1s_{n+1-i} = u$.
\item[\emph{(iii)}] If $w = u_1s_{n+1-i}u_2$ for some $i \in [2,k-1]$,
where $u_1 \suff s_{n+1-i}$ and  $u_2 \pref s_n$, then $u_1 = u$
or $u_1s_{n+1-i} = u$.
\end{itemize}
\end{lemma}
\begin{proof} 
(i) Other than the case when $u_1 = u$, $u_2 = s_n$ and $v =
\empt$, the only other possibility is: \vspace{0.4cm}

\centerline{\includegraphics{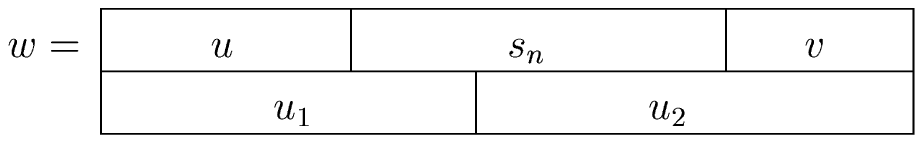}}

(Note that $u_1 \suff s_{n+1-i}$ for some $i \in [1,k]$, and
therefore $|u_1| \leq |s_{n+1-i}| \leq |s_n|$.)

In this case, using the figure, we write $u_1 = uu'$, $s_n = u'v'$, $u_2 = v'v$ for some $u'$, $v'$ ($u' \ne \empt$). As $v'$ is a prefix of $s_n$, we have $s_n = v'v''$ for some $v''$, thus $u'$ and $v''$ are conjugate. So there exist $e$, $f$ and non-negative integers $p$, $q$ such that $v' = (ef)^pe$, $u' = (ef)^q$, and 
$v'' = (fe)^q$ with $ef$ primitive. Hence $s_n = (ef)^{p+q}e$. As $u'$ is a suffix of $s_{n+1-i}$ which is a prefix of $s_n$, we must have, by primitivity of $ef$, $s_{n+1-i} = (ef)^r$, and then $r = 1$. But $u'$ is non-empty, so $u' = ef = s_{n+1-i}$, and it follows that $u = \empt$; a contradiction.

(ii) Let $i \in [2,k]$ be a fixed integer. Since $u_1 \suff
s_{n+2-i}$ and $u_2 \pref s_n$, we have $|u_1| \leq |s_{n+2-i}|
\leq |s_n|$ and $|u_2| \leq |s_n|$. Accordingly, there exist only
three possibilities (other than $u_1 = u$ or $u_1s_{n+1-i} = u$),
and these are: \vspace{0.7cm}

\centerline{\includegraphics{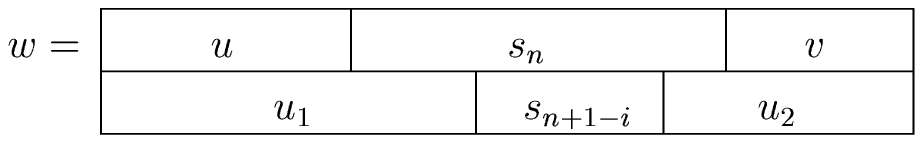}}

or

\centerline{\includegraphics{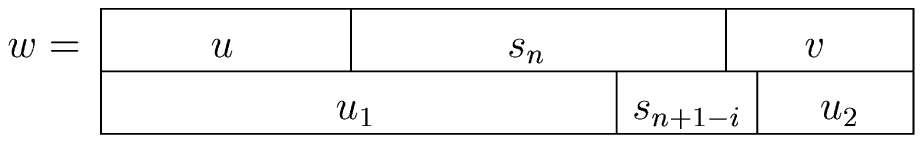}}

or

\centerline{\includegraphics{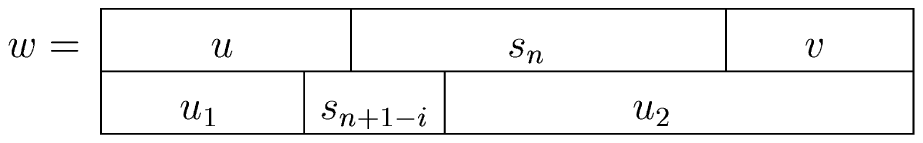}}

In the first instance, $us_n v = u_1s_{n+1-i}u_2 = uu^\prime
s_{n+1-i}v^\prime v$, where $u_2 = v^\prime v$ with $v^\prime
\suff s_n$ and $u_1 = uu^\prime$ with $u^\prime \pref s_n$. That
is, $s_n = u^\prime s_{n+1-i}v^\prime$, where $u^\prime \pref s_n$
and $v^\prime \suff s_n$, and $u_1 = uu^\prime \suff s_{n+2-i}$.
Therefore, $u^\prime \suff s_{n+2-i}$, and hence the word
$s_{n+1-i}$ must be preceded by the last letter of $s_{n+2-i}$.
However, since $u^\prime$ is also a prefix of $s_n =
s_{n-1}^{d_{n}}\cdots s_{n+1-k}^{d_{n+2-k}}s_{n-k}$, where
$s_{n-1}$, $\ldots$~, $s_{n+1-k}$, $s_{n-k}$ do not share a common
last letter (by Proposition \ref{P:27.11.03(1)}), one is forced to
presume that $u^\prime = s_{n-1}^{d_{n}}s_{n-2}^{d_{n-1}}\cdots
s_{n+2-i}^{d_{n+3-i}}$ (resp.~$u^\prime = \empt$) when $i \in
[3,k]$ (resp.~$i = 2$). This contradicts the fact that $1\leq
|u^\prime| < |s_{n+2-i}|$.

In the second instance, we have $us_nv = u_1 s_{n+1-i}u_2 = uu^\prime
s_{n+1-i}u_2$, where $u_1 = uu^\prime$ with $u^\prime \pref s_n$
and $u_2 \pref s_n$. Consider the word $w^\prime := w{u_2}^{-1} =
us_nv{u_2}^{-1} = us_nv^\prime$, i.e., $w^\prime = us_nv^\prime =
u_1s_{n+1-i}$, where $v^\prime \pref v$ and $v^\prime \suff
s_{n+1-i}$. Since $u_1 \suff s_{n+2-i}$ and $s_{n+1-i} \pref s_n$,
it follows from assertion (i) that $u_1 = u$ and hence $s_{n+1-i}
= s_nv^\prime$, which is absurd unless $i = 1$ and $v^\prime =
\empt$. But $i > 1$, so this situation is impossible.

Lastly, $us_n v = u_1s_{n+1-i}u_2 = u_1s_{n+1-i}v^\prime v$, where
$u_2 = v^\prime v \pref s_n$ with $v^\prime \suff s_n$. Consider
the word $w^\prime := u_1^{-1}w = u_1^{-1}us_nv = u^\prime s_n v$,
i.e., $w^\prime = u^\prime s_n v = s_{n+1-i}u_2$, where $u^\prime
\suff u$ and $u^\prime \pref s_{n+1-i}$. Since $u_2 \pref s_n$,
one obtains, as an immediate consequence of claim (i), $u^\prime =
s_{n+1-i}$, $u_2 = s_{n}$, $v = \empt$, and hence $u =
u_1s_{n+1-i}$; a contradiction since $|u_1s_{n+1-i}| > |u_1|$.

One can prove assertion (iii) in a similar manner.
\end{proof}

\begin{lemma} \label{L:26.07.04(2)}
Let $c \in \cA_k$ and $n \in \NN$ be fixed. Consider an occurrence
of $cs_n$ in $\bs$. Then the letter $c$ is the last letter of a
block $s_{n+1-i}$ of the $n$-partition of $\bs$, for some $i \in
[1,k]$, and the integer $i$ $($equiv. the block $s_{n+1-i})$ is
uniquely determined by $c$. In particular, in every occurrence of
$s_{n+1-i}s_n$ in $\bs$, the word $s_{n+1-i}$ is a block in the
$n$-partition of $\bs$.
\end{lemma}
That is, occurrences of words $w$ containing $cs_n$ ($c \in
\cA_k$) must be aligned to the $n$-partition of $\bs$.

\begin{proof} 
This assertion follows from Lemma \ref{L:26.07.04(1)}. The case
$n=0$ is trivial, and for $n\geq 1$, observe from Lemma
\ref{L:26.07.04(1)} that the given $s_n$ is either an $s_n$ block
in the $n$-partition of $\bs$ or has an $s_{n+1-j}$ block of the
$n$-partition as a prefix, for some $j \in [2,k]$. In the first
case, to the left of $s_n$ there is an $s_{n+1-l}$ block, for some
$l \in [1,k]$. Whereas, in the second case, there is an
$s_{n+1-j}$ or $s_{n+2-j}$ block (of the $n$-partition) to the
left of $s_n$. That is, $s_n$ is preceded by an $s_{n+1-i}$ block
of the $n$-partition for some $i \in [1,k]$. Since the last
letters of $s_{n}$, $s_{n-1}$, $\ldots$~, $s_{n+1-k}$ are mutually distinct
(by Proposition \ref{P:27.11.03(1)}), it is clear that $i$ (and
hence $s_{n+1-i}$) is uniquely determined by the letter $c$.
\end{proof}

We can now determine the exact index of $s_n$ in $\bs$.

\begin{lemma} \label{L:06.12.03(2)}
For any $n \geq 1$, the word of maximal length that is a factor of
both $\bs$ and the infinite sequence $(s_n)^\omega :=
s_ns_ns_n\cdots$ is $s_n^{d_{n+1}+2}D_{n-k}$, i.e.,
\emph{ind}$(s_n) = 2 + d_{n+1} + |D_{n-k}|/|s_n|$.
\end{lemma}
\begin{proof}  
According to Lemma \ref{L:26.07.04(2)}, any occurrence of $s_n^p$
($p\geq 2$) must be aligned to the $n$-partition of $\bs$. By
inspection of the $n$-partition of $\bs$ (see
\eqref{eq:n-part-k}), it is not hard to see that between two
successive $s_{n+1-k}$ blocks there is a word possessing one of
the following $k$ forms:
\[ 
s_n^qs_{n}^{d_{n+1}}s_{n-1}^{d_n}\cdots s_{n+2-k}^{d_{n+3-k}},
\quad q \in \{0,1\},
\]
or
\[
s_{n}^{d_{n+1}}s_{n-1}^{d_n}\cdots s_{n+2-i}^{d_{n+3-i}}
s_{n+1-i}s_{n}^{d_{n+1}}s_{n-1}^{d_n}\cdots s_{n+2-k}^{d_{n+3-k}},
\quad i \in [2,k-1].
\]
Thus, the alignment property implies that an occurrence of $s_n^p$ ($p\geq 2$) is
either a prefix of
\begin{equation} \label{eq:28.07(1)}
s_{n}^{d_{n+1}+r}s_{n-1}^{d_n}\cdots
s_{n+2-k}^{d_{n+3-k}}s_{n+1-k}z_1
\end{equation}
for some integer $r \leq 1$ and suitable $z_1$,  or a prefix of
\begin{equation} \label{eq:28.07(2)}
s_{n}^{r}s_{n-1}^{d_n}\cdots s_{n+2-i}^{d_{n+3-i}}
s_{n+1-i}s_{n}^{d_{n+1}}s_{n-1}^{d_n}\cdots s_{n+2-k}^{d_{n+3-k}}
s_{n+1-k}z_2
\end{equation}
for some $i \in [2,k-1]$, $r \leq d_{n+1}$ and suitable $z_2$.

Now, suppose $s_n^p$ is a prefix of the word \eqref{eq:28.07(1)}.
Since $s_{n-1}^{d_n}\cdots s_{n+2-k}^{d_{n+3-k}}s_{n+1-k}s_n$ is
not a prefix of $s_ns_n$ (in fact, it is the word
$s_n(s_{n+1-k}^{d_{n+2-k}-1}s_{n-k})^{-1}s_n = s_nt_{n-1}$),
$s_n^p$ must also be a prefix of
\begin{equation} \label{eq:28.07(3)}
  s_{n}^{d_{n+1}+r}s_{n-1}^{d_n}\cdots s_{n+2-k}^{d_{n+3-k}}s_{n+1-k}s_n
  = s_{n}^{d_{n+1}+r+1}t_{n-1}.
\end{equation}
As in the proof of Lemma \ref{L:index1}, one can show that the
prefix index of $s_n$ in the word \eqref{eq:28.07(3)} is $d_{n+1}
+ r + 1 + |D_{n-k}|/|s_n|$, which is at most $d_{n+1} + 2 +
|D_{n-k}|/|s_n|$. Furthermore, in the word \eqref{eq:28.07(2)}, it
is clear that the prefix index of $s_n$ is less than for
\eqref{eq:28.07(1)} (since $s_{n-1}^{d_n}\cdots
s_{n+2-i}^{d_{n+3-i}}s_{n+1-i}s_{n}$ has length less than the word
$s_{n-1}^{d_n}\cdots s_{n+2-k}^{d_{n+3-k}}s_{n+1-k}s_n$ and is not
a prefix of $s_ns_n$). Whence, it has been shown that ind$(s_n)
\leq d_{n+1} + 2 + |D_{n-k}|/|s_n|$, and so the result is now 
an easy consequence of Lemma \ref{L:29.06.04(2)} (which gives
ind$(s_n) \geq d_{n+1} + 2 + |D_{n-k}|/|s_n|$).
\end{proof}

The following analogue of Lemma 3.5 in \cite{dDdL02thei} is
required in order to prove Theorem \ref{T:22.11.03(5)}.

\begin{lemma} \label{L:18.05.04(3)} Let $n \in \NN^+$ 
and suppose $u \prec \bs$ with $|s_{n}| \leq |u| < |s_{n+1}|$.
Then the following assertions hold.
\begin{itemize}
\item[$(1)$] For all $i \in [1, k]$, if $u$ starts at position $l$ in some $s_{n+1-i}$ block in the $n$-partition of $\bs$ and also starts at position $m$ in some factor $s_{n+1-i}$ of $\bs$, then $l = m$.  
\item[$(2)$] For all $i \in [1,k-1]$, if $u$ can start at position $l$ in $s_{n+1-i}$ and at position
$m$ in $s_{n-i}$, 
then $l = m$.
\end{itemize}
\end{lemma}
\begin{proof} 
By inspection of the $n$-partition of $\bs$, notice that, for $1
\leq i \leq k-1$, an $s_{n+1-i}$ block is followed by either an
$s_{n+1-i}$ block, an $s_{n}$ block, or an $s_{n-i}$ block.
Furthermore, an $s_{n+1-k}$ block is always followed by an $s_n$
block.

Let $u_{n+1-i}$ be the prefix of $u$ of length $|s_{n+1-i}|$.

(1) Let $1 \leq i \leq k$ and consider an occurrence of $u$
that starts in an $s_{n+1-i}$ block of the $n$-partition of $\bs$.
If this $s_{n+1-i}$ block is followed by an $s_{n+1-i}$ block,
then $u_{n+1-i}$ is a conjugate of $s_{n+1-i}$ as $u_{n+1-i} \prec
s_{n+1-i}s_{n+1-i}$ and $|u_{n+1-i}| = |s_{n+1-i}|$. Similarly, if
this $s_{n+1-i}$ block is followed by an $s_{n}$ block, then
$u_{n+1-i}$ is a conjugate of $s_{n+1-i}$ since $s_{n+1-i} \pref
s_n$. And, if this $s_{n+1-i}$ block is followed by an $s_{n-i}$
block, then again $u_{n+1-i}$ is a conjugate of $s_{n+1-i}$.
Indeed, in the $(n+k-i)$-partition of $\bs$, $s_{n+1-i}$ is always
followed by an $s_{n+k-i}$ block, which has $s_{n+1-i}$ as a
prefix; whence $u_{n+1-i} \prec s_{n+1-i}s_{n+1-i}$. So, in any
case, $u_{n+1-i}$ is a conjugate of $s_{n+1-i}$, and the result
follows from the fact that the conjugates of $s_{n+1-i}$ are
distinct. 

(2) Let $1 \leq i \leq k-1$. Suppose the word $u$ has occurrences
starting in $s_{n+1-i}$ blocks as well as $s_{n-i}$ blocks in the
$n$-partition of $\bs$. (Note that this implies $n\geq i$.) First
consider an occurrence of $u$ beginning in a block of the form
$s_{n-i}$ of the $n$-partition. As an $s_{n-i}$ block is always
followed by an $s_{n+k-i-1}$ block in the $(n+k-i-1)$-partition of
$\bs$ and $s_{n+1-i} \pref s_{n+k-i-1}$, we have
\[
  u_{n+1-i} \prec s_{n-i}s_{n+1-i} = s_{n-i}s_{n-i}^{d_{n-i+1}}\cdots
s_{n+2-i-k}^{d_{n+3-i-k}}s_{n+1-i-k}.
\]
Thus, in light of Lemma \ref{L:26.07.04(2)}, we have the following
fact:
\begin{equation} \label{eq:30.07(1)}
cs_{n-i} \not\prec u_{n+1-i} \quad \mbox{where $c \in \cA_k$ and
$c \suff s_{n+1-i-k}$}.
\end{equation}

Consider an occurrence of $u$ starting in an $s_{n+1-i}$ block
of the $n$-partition, which can be factorized as
\begin{equation} \label{eq:30.07(2)}
  s_{n+1-i} = s_{n-i}^{d_{n-i+1}}s_{n-i-1}^{d_{n-i}}\cdots
  s_{n+2-i-k}^{d_{n+3-i-k}}s_{n+1-i-k}.
\end{equation}
We distinguish two cases, below.

\emph{Case $1$}: The word $u$ begins in the left-most $s_{n-i}$
block in \eqref{eq:30.07(2)} when $d_{n-i+1}\geq 2$. In this case,
$u_{n-i}$ is a conjugate of $s_{n-i}$ and hence, as deduced in
(1), the starting position of $u$ in this $s_{n-i}$ block must  
coincide with
its starting position in any occurrence of $s_{n-i}$ in the
$n$-partition of $\bs$.

\emph{Case $2$}: The word $u$ does not start in the left-most
$s_{n-i}$ block in \eqref{eq:30.07(2)}. The block to the right of
$s_{n+1-i}$ in the $n$-partition is either another $s_{n+1-i}$, or
an $s_{n-i}$, or an $s_n$. In any case, $s_{n-i}$ is a prefix of
this block to the right of $s_{n+1-i}$, which implies $u_{n+1-i}$
contains an occurrence of $s_{n+1-i-k}s_{n-i}$. This  contradicts 
\eqref{eq:30.07(1)}.
\end{proof}

\begin{proof}[Proof of Theorem $\ref{T:22.11.03(5)}$] 
Clearly, $p(m;l_1) \geq p(m;l_2)$ if $l_1 \leq l_2$. Thus, it
suffices to show that for $m \not\in \bigcup_{i=1}^k \cD_i(n)$, we
have
\begin{equation} \label{eq:pm2-k}
p(m;2) = 0,
\end{equation}
i.e., there are no squares of words of length $m$ in $\bs$.

Suppose \eqref{eq:pm2-k} does not hold for some $m$ satisfying
\begin{equation} \label{eq:D1UD2}
  m \not\in \bigcup_{i=1}^k \cD_i(n),
\end{equation}
and let $u$ be a word of length $m$ with $|s_n| \leq m <
|s_{n+1}|$ such that $u^2 \prec \bs$. For convenience, we shall
write $u^2 = u^{(1)}u^{(2)}$  
to allow us to refer to the two separate occurrences of $u$. Let $1 \leq i,j \leq
k$. Obviously, $u^{(1)}$ starts at position $q$, say, 
in some $s_{n+1-i}$ block of the $n$-partition of $\bs$. Further, by Lemma
\ref{L:18.05.04(3)}, $u^{(2)}$ also starts in some $s_{n+1-j}$
block of the $n$-partition of $\bs$ at position $q$. From the
proof of Lemma \ref{L:06.12.03(2)}, recall that two consecutive
$s_{n+1-k}$ blocks in the $n$-partition of $\bs$ are separated by
a word of one of the following $k$ forms:
\[ 
s_n^rs_{n}^{d_{n+1}}s_{n-1}^{d_n}\cdots s_{n+2-k}^{d_{n+3-k}},
\quad r \in \{0,1\},
\]
or
\[
s_{n}^{d_{n+1}}s_{n-1}^{d_n}\cdots s_{n+2-i}^{d_{n+3-i}}
s_{n+1-i}s_{n}^{d_{n+1}}s_{n-1}^{d_n}\cdots s_{n+2-k}^{d_{n+3-k}},
\quad i \in [2,k-1].
\]
If we also keep in mind that $|s_n| \leq |u| < |s_{n+1}|$, then
using Lemma \ref{L:18.05.04(3)} we see that the possible lengths
$|u|$ of $u$ are:
\[
  |s_n^r| \quad \mbox{and} \quad
  |s_n^rs_{n-1}^{d_{n}}\cdots s_{n-i+1}^{d_{n-i+2}}s_{n-i}| 
\]
where $1 \leq i \leq k-1$ and $1 \leq r \leq d_{n+1}$ (with $r \ne
d_{n+1}$ if $i = k-1$ as $|u| < |s_{n+1}|$). 
Therefore, $m$ does not satisfy \eqref{eq:D1UD2}; a contradiction.
\end{proof}

The next five propositions, which have some interest in themselves, are needed in the next two sections where we shall prove our main results concerning squares, cubes, and higher powers in $\bs$. 

\begin{notation}
Given $l \in \NN$ and $w \in \cAstar_k$, denote by Pref$_l(w)$ the
prefix of $w$ of length $l$ if $|w| \geq l$, $w$ otherwise.
Likewise, denote by Suff$_l(w)$ the suffix of $w$ of length $l$ if
$|w|\geq l$, $w$ otherwise.

Recall that $\Omega_n^r$ denotes the set of singular $n$-words of the 
$r$-th kind ($1\leq r\leq k-1$), as defined in Theorem \ref{T:06.09.04(3)}.
\end{notation}

\begin{proposition} \label{P:19.05.04(1)} Let $n \in \NN^+$.   
Suppose $w \in \Omega_{n+1-i}^1$ for some $i \in [1, k-1]$ and 
let $v =$ \emph{Pref}$_l(w)$ where $1 \leq l \leq |G_{n+1-i,1}| - 1$.
Then the word $vs_{n+1-i}$ occurs at position $p$ in $\bs$ if and only if the
$n$-partition of $\bs$ contains an $s_{n}$ starting at position
$p+l$ and an $s_{n-i}$ ending at position $p+l-1$. In particular,
$w$ occurs at exactly those positions where $vs_{n+1-i}$ occurs in
$\bs$.
\end{proposition}
\begin{proof}  
Let $i \in [1, k-1]$ be fixed and let $1 \leq l \leq |G_{n+1-i,1}| - 1$.

First note that $|w| = |s_{n+1-i}|$ and $w \prec
x^{-1}\rev{G}_{n+1-i,1}D_{n-i}G_{n+1-i,1}x^{-1}$ where $x \in \cA_k$,  
by definition of $\Omega_{n+1-i}^1$.  Since $|D_{n-i}G_{n+1-i,1}| = |s_{n+1-i}|$, the word $v = \mbox{Pref}_l(w)$ is a suffix of  $x^{-1}\rev{G}_{n+1-i,1}$ which, in turn, is a suffix of $s_{n-i}$ as $\rev{G}_{n+1-i,1} =  G_{n-i,k-1}$. 

Now, by Lemma \ref{L:18.05.04(3)}, the word $s_{n+1-i}$ can only
occur at the starting positions of blocks (in the $n$-partition)
of the form $s_n$, $s_{n-1}$, $\ldots$~, $s_{n+1-k}$, all of which
have different last letters (by Proposition \ref{P:27.11.03(1)}).
In particular, each $s_{n-j}$ block ($0\leq j\leq k-1$, $j\ne i$)
of the $n$-partition of $\bs$ has a different last letter to
$s_{n-i}$ (and hence $v$). One should note, however, that an
$s_{n-i}$ block of the $n$-partition is never followed by an
$s_{n+1-i}$ block (except if $i =1$, in which case we do 
have certain $s_{n}$ blocks preceded by $s_{n-1}$ blocks). Also 
observe that if $z = \mbox{Suff}_l(s_{n+1-i})$, then an $s_{n-i}$ block of the $n$-partition is only ever
followed by $s_{n+1-i}z^{-1} =D_{n-i}G_{n+1-i,1}z^{-1}$ if it
is followed by an $s_n$ block of the $n$-partition.   
Taking all of this into account, one deduces that the word $vs_{n+1-i}$ occurs only at positions in $\bs$
where an $s_{n-i}$ block of the $n$-partition is followed by an
$s_n$ block, which has $s_{n+1-i}$ as a prefix. This completes the proof of the 
first assertion.

As for the second assertion, recall that $w$ begins with the word $v$ which is a 
non-empty suffix of $x^{-1}\rev{G}_{n+1-i,1} = x^{-1}G_{n-i,k-1} \suff s_{n-i}$. Consequently, $w$ occurs at every $(|s_{n-i}| - l + 1)$-position of an $s_{n-i}$ block that is
followed by an $s_n$ block in the $n$-partition of $\bs$, i.e., $w$ occurs 
where the prefix $vs_{n+1-i}$ of $vs_n$ occurs in $\bs$. By Lemma
\ref{L:18.05.04(3)}, the only other position where $w$ may occur
(besides where $vs_{n+1-i}$ occurs) is in the
$(|s_{n-i}|-l+1)$-position of an $s_{n}$ block that is preceded by 
an $s_{n-i}$ block. Now, to the right of this type of $s_{n}$
block (in the $n$-partition) there appears another $s_{n}$ block
or an $s_{n-1}$ block. The fact that $s_{n+1-i}s_{n-i} \pref s_ns_{n-1} \pref s_ns_n$ 
implies that $w$ ends with the prefix of $s_{n-i}$ of length $|s_{n-i}| - l$. More precisely, $w$
ends with the word
\[ 
  D_{n+1-i-k}z_1 
\]  
where $z_1$ is a non-empty prefix of $\rev{G}_{n+1-i,1}$ of length
$|z_1| = |G_{n+1-i,1}| - l$. On the other hand, by definition of $w$, we have that
$w$ ends with 
\[
  D_{n-i}z_2 
\]
where $z_2$ is a non-empty prefix of $G_{n+1-i,1}$ of length $|z_2| = |G_{n+1-i,1}| - l$. It is impossible for
both situations to occur, so we conclude that $w$ occurs at
exactly those positions where $vs_{n+1-i}$ occurs.
\end{proof}

\begin{notation}
For $n\geq 1$, denote by $\PP_n$ the set of all formal
positions of $s_n^{d_{n+1}-1}s_{n-1}$ in the
$(n-1)$-partition of $\bs$.
\end{notation}

\begin{proposition} \label{P:25.05.04(1)} For any $n\in \NN^+$,
the set of all positions of $D_{n}$ in $\bs$ is $\PP_n$.
\end{proposition}
\begin{proof} 
Proof is by induction on $n$. For $n = 1$, $D_n = D_1 =
s_1^{d_2-1}s_0^{d_1}$, and hence $D_1$ occurs at exactly those
places in $\bs$ where $s_1^{d_2-1}s_{0} =
(a_1^{d_{1}}a_2)^{d_{2}-1}a_1$ occurs in the $0$-partition of
$\bs$. We claim that there is a one-to-one correspondence from the
set of all positions of $D_n$ in the $(n-1)$-partition of $\bs$
to the set of all positions of $D_{n+1}$ in the $n$-partition of
$\bs$ (see \eqref{eq:n-part-k}). Assume that $\PP_n$ gives all of
the occurrences of $D_n$ in the $(n-1)$-partition of $\bs$. Since
$D_{n+1} = s_{n+1}^{d_{n+2}-1}s_{n}D_n =
D_n\rev{s}_n(\rev{s}_{n+1})^{d_{n+2}-1}$, $D_{n+1}$ occurs at any
place in $\PP_{n+1}$. Conversely, since each occurrence of
$D_{n+1}$ in \eqref{eq:n-part-k} naturally gives rise to an
occurrence of $D_n$ in the $(n-1)$-partition of $\bs$, the word
$D_{n+1}$ must occur in $\bs$ at exactly those places given by
$\PP_{n+1}$.
\end{proof}

Consider two distinct occurrences of a factor $w$ in $\bs$, say 
\[
  \bs = uw\mathbf{v} = u^\prime w\mathbf{v}^\prime, \quad |u^\prime| > |u|,
\]
where $\mathbf{v}$, $\mathbf{v}^\prime \in \cA_k^\omega$. These two 
occurrences of $w$ in $\bs$ are said to be \emph{positively separated} 
(or \emph{disjoint}) if $|u^\prime| > |uw|$, in which case 
$u^\prime = uwz$ for some $z \in \cA_k^+$, and hence 
$\bs = uwzw\mathbf{v}^\prime$.

\begin{proposition} \label{P:06.12.03}
For any $n \in \NN^+$, successive occurrences of a singular word
$w \in \bigcup_{j=1}^{k-1}\Omega_n^j$ in $\bs$ are positively
separated.
\end{proposition}
\begin{proof} 
Let $1\leq r\leq k-1$. For $1 \leq n \leq r$, observe that
\[
  \Omega_n^r = \{w \in \cA_k^* ~: ~|w| = |s_n| ~\mbox{and} ~w \prec 
  s_{n-1}D_{n-1}D_{n-r}^{-1}s_{n-1}D_{n-1} \},
\]
where $D_{n-r}^{-1} = a_{n-r+k+1}$ for $1\leq n < r$. 
It is left to the reader to verify that consecutive occurrences of
a word $w \in \Omega_n^r$ ($1\leq n \leq r$) are positively separated in $\bs$.

Now take $n \geq r+1$ and suppose $w \in \Omega_n^r$. Then $D_{n-r}$
will occur in $w$. By Propositions \ref{P:25.05.04(1)}, the word
$D_{n-r}$ occurs at exactly those places where 
$s_{n-r}^{d_{n-r+1}-1}s_{n-r-1}$ occurs in the $(n-r-1)$-partition
of $\bs$. 
First note that the letter just before $D_{n-r}$ in $w$ is the
last letter of $\rev{G}_{n,r}$, which is the first letter
$G_{n,r}$, and hence the last letter of $s_{n-r-k}$ (by
Propositions \ref{P:27.11.03(1)} and \ref{P:19.04.04(2)}). On the
other hand, in the word $w$, the letter just after $D_{n-r}$ is
the first letter of $G_{n,r}$. Since there are $k$ different
return words of $s_{n-r}^{d_{n-r+1}}s_{n-r-1}$ in $\bs$, there
exist $k$ different possibilities for occurrences of
$s_{n-r}^{d_{n-r+1}}s_{n-r-1}$ in the $(n-r-1)$-partition of
$\bs$; namely:
\begin{itemize}
\item[(1)]

\begin{align*}
&(s_{n-r}^{d_{n-r+1}-1}s_{n-r-1})s_{n-r-1}^{d_{n-r}-1}
  s_{n-r-2}^{d_{n-r-1}}\cdots s_{n-r-k+1}^{d_{n-r-k+2}}s_{n-r-k}
  (s_{n-r}^{d_{n-r+1}-1}s_{n-r-1}) \\ &= ~D_{n-r}D_{n-r-k}^{-1}
  D_{n-r}D_{n-r-1}^{-1} \\
  &= ~D_{n-r}(\rev{s}_{n-r})^{d_{n-r+1}}D_{n-r-1}^{-1} \\
  &= ~D_{n-r}(\rev{s}_{n-r})^{d_{n-r+1}-1}\rev{G}_{n-r,1};
\end{align*}
\item[(2)] 
\begin{align*}
&(s_{n-r}^{d_{n-r+1}-1}s_{n-r-1})s_{n-r-1}^{d_{n-r}-1}
  s_{n-r-2}^{d_{n-r-1}}\cdots s_{n-r-l+1}^{d_{n-r-l+2}}s_{n-r-l}
  (s_{n-r}^{d_{n-r+1}-1}s_{n-r-1}) \\ &= ~D_{n-r}D_{n-r-l}^{-1}
  (s_{n-r}^{d_{n-r+1}-1}s_{n-r-1}) \\
  &= ~\begin{cases}
      D_{n-r}G_{n-r,l}{s}_{n-r}^{d_{n-r+1}-2}s_{n-r-1} 
      &\mbox{if $d_{n-r+1} \geq 2$}, \\
      D_{n-r}G_{n-r-1,l-1} &\mbox{if $d_{n-r+1} = 1$}, 
     \end{cases}
\end{align*}
where $2 \leq l \leq k-1$; 
\item[(3)]
\begin{align*}
&(s_{n-r}^{d_{n-r+1}-1}s_{n-r-1})(s_{n-r}^{d_{n-r+1}-1}s_{n-r-1})
s_{n-r-1}^{d_{n-r}-1}
  s_{n-r-2}^{d_{n-r-1}}\cdots s_{n-r-k+2}^{d_{n-r-k+3}}s_{n-r-k+1} \\
  &= ~\begin{cases}
      D_{n-r}\rev{s}_{n-r-1}(\rev{s}_{n-r})^{d_{n-r+1}-2}\rev{G}_{n-r,k-1}
      &\mbox{if $d_{n-r+1} \geq 2$}, \\
      D_{n-r}\rev{G}_{n-r-1,k-2} ~~~ \qquad
      &\mbox{if $d_{n-r+1} = 1$}.
     \end{cases}
\end{align*}
\end{itemize}
Thus, if $d_{n-r+1} \geq 2$, the word $D_{n-r}$ is followed by
either $\rev{s}_{n-r}$, $G_{n-r,l}$, or $\rev{s}_{n-r-1}$, of
which only $\rev{s}_{n-r}$ has the same first letter as $G_{n,r}$.
Similarly, if $d_{n-r+1} = 1$, the word $D_{n-r}$ is followed by
either $\rev{G}_{n-r,1}$, $G_{n-r-1,l-1}$, or
$\rev{G}_{n-r-1,k-2}$, of which only $\rev{G}_{n-r,1}$ has the
same first letter as $G_{n,r}$. Therefore, only in case (1) will
we have $D_{n-r}$ followed by the first letter of $G_{n,r}$.
Accordingly, one deduces that any occurrence of $w$ in $\bs$
corresponds to a formal occurrence of the word
\[
  s_{n-r-k}(s_{n-r}^{d_{n-r+1}-1}s_{n-r-1})s_{n-r-1}^{d_{n-r}-1}
  s_{n-r-2}^{d_{n-r-1}}\cdots s_{n-r-k+1}^{d_{n-r-k+2}}s_{n-r-k}
  (s_{n-r}^{d_{n-r+1}-1}s_{n-r-1})
\]
in the $(n-r-1)$-partition of $\bs$. Hence, we conclude that
occurrences of $w$ are positively separated in $\bs$ since a word
of the above form is positively separated in the
$(n-r-1)$-partition.
\end{proof}

The next proposition follows from Lemma \ref{L:06.12.03(2)}, and Propositions  
\ref{P:19.05.04(1)} and \ref{P:06.12.03}.

\begin{proposition} \label{P:27.11.03(2)} 
Let $n \in \NN^+$ and suppose $u \prec \bs$ with $|u| = |s_n|$.
Then $u^2 \prec \bs$ if and only if $u \in \cC(s_n)$. In particular,
if $u$ is a singular word of any kind of $\bs$, then $u^2 \nprec
\bs$. Moreover, for any $n\geq k-1$, if $u^2 \prec \bs$ with
$|s_n| \leq |u| < |s_{n+1}|$, then $u$ does not contain a singular
word from the set ~$\Omega_{n+2-k}^1$.
\end{proposition}
\begin{proof} 
As $s_n^{d_{n+1}+2}D_{n-k} \prec \bs$ (see Lemma
\ref{L:06.12.03(2)}), the square of any conjugate of $s_n$ is a
factor of $\bs$. (Note that $s_n^{d_{n+1}+2}D_{n-k} =
s_n^{d_{n+1}+2}a_{n+1}^{-1}$ for $1\leq n\leq k-1$.) Now recall
that the set of all factors of $\bs$ of length $|s_n|$ is the
disjoint union of the sets $\cC(s_n)$ and
$\bigcup_{j=1}^{k-1}\Omega_{n}^j$. Consequently, the first two
assertions are deduced from Proposition \ref{P:06.12.03}.

For the last statement, let $n\geq k-1$ and suppose $u^2 =
u^{(1)}u^{(2)}$ is an occurrence of $u^2$ in $\bs$, where $|s_n|
\leq |u| < |s_{n+1}|$. Also assume $w \in \Omega_{n+2-k}^1$ and $w
\prec u$. Clearly, $w$ occurs in both $u^{(1)}$ and $u^{(2)}$ at
the same position. By Proposition \ref{P:19.05.04(1)} (with $i = k-1$),
different occurrences of $w$ correspond to different occurrences
of $s_{n+1-k}$ blocks in the $n$-partition of $\bs$ (as an
$s_{n+1-k}$ block is always followed by an $s_n$ block). Between
two consecutive $s_{n+1-k}$ blocks in the $n$-partition, there is
a word taking one of the following $k$ forms:
\[ 
s_n^rs_{n}^{d_{n+1}}s_{n-1}^{d_n}\cdots s_{n+2-k}^{d_{n+3-k}},
\quad r \in \{0,1\},
\]
or
\[
s_{n}^{d_{n+1}}s_{n-1}^{d_n}\cdots s_{n+2-i}^{d_{n+3-i}}
s_{n+1-i}s_{n}^{d_{n+1}}s_{n-1}^{d_n}\cdots s_{n+2-k}^{d_{n+3-k}},
\quad i \in [2,k-1].
\]
Therefore, the distance between consecutive occurrences of
$s_{n+1-k}$ blocks in the $n$-partition of $\bs$ is
\[
  |s_{n+1}| ~~(r = 0), \quad |s_{n+1}| + |s_n| ~~(r=1), \quad \mbox{or}
  \quad |s_{n+1}| - (|D_{n+1-i}| - |D_{n+1-k}|) + |s_{n+1}|,
\]
and all of these distances are at least $|s_{n+1}|$, which implies
$|u| \geq |s_{n+1}|$; a contradiction.
\end{proof}

More generally, we have the following proposition.

\begin{proposition} \label{P:04.08.04(1)}  
Let $n \in \NN^+$  and suppose $u^2 \prec \bs$
with $|s_n| \leq |u| < |s_{n+1}|$. Then $u$ does not contain a
singular word from the set $\Omega_{n+1-i}^1$ for any $i \in
[1,k-1]$.
\end{proposition}
\begin{proof}  
The case when $i = k-1$ is proved in Proposition \ref{P:27.11.03(2)}, 
so take $i \in [1,k-2]$. Let $u^2 = u^{(1)}u^{(2)}$ be an
occurrence of $u^2$ in $\bs$, where $|s_n| \leq |u| < |s_{n+1}|$.
Assume $w \in \Omega_{n+1-i}^1$ for some $i \in [1,k-2]$, and $w
\prec u$. Clearly, $w$ occurs in both $u^{(1)}$ and $u^{(2)}$ at
the same position. By Proposition \ref{P:19.05.04(1)}, different
occurrences of $w$ correspond to different occurrences of
$s_{n-i}$ blocks that are followed by $s_n$ blocks in the
$n$-partition of $\bs$. By inspection of the $n$-partition (see
\eqref{eq:n-part-k}), the word of minimal length that separates
two such $s_{n-i}$ blocks is
\[
  s_n^{d_{n+1}}s_{n-1}^{d_n}\cdots s_{n+2-k}^{d_{n+3-k}}s_{n+1-k}
  s_n^{d_{n+1}}s_{n-1}^{d_n}\cdots s_{n+1-i}^{d_{n+2-i}}.
\]
That is, the minimal distance between two consecutive occurrences
of an $s_{n-i}$ block (with each appearance followed by an $s_n$
block) is 
$|s_{n+1}|+|s_n| + |D_n| - |D_{n-i}| > |s_{n+1}| + |s_n|$,  
which implies $|u| > |s_{n+1}| + |s_n|$; a contradiction.
\end{proof}

\subsection{Squares} \label{SS:squares_in_s}

The next two main theorems 
concern squares of factors of $\bs$ of length
$m < d_1+1 = |s_1|$ and length $m \geq |s_1|$, respectively.

A letter $a$ in a finite or infinite word $w$ is said to be 
\emph{separating for $w$} if any factor of length 2 of $w$
contains the letter $a$. For example, $a$ is separating for the infinite word $(aaba)^\omega$.  
If $a$ is separating for an infinite word
$\bx$, then it is clearly separating for any factor of $\bx$. According to 
\cite[Lemma 4]{xDjJgP01epis}, since the standard episturmian word
$\bs$ begins with $a_1$, the letter $a_1$ is separating for $\bs$
and its factors.

\begin{theorem} \label{T:02.08.04(1)}
For $1 \leq r \leq d_1$, we have
\[
  p(r;2) = \begin{cases}
            1 \quad \mbox{if $r \leq (d_1+1)/2$}, \\
            0 \quad \mbox{if $r > (d_1+1)/2$}.
           \end{cases}
\]
In particular, $\cP(r;2) = \{(a_1^r)^2\}$ for $r \leq (d_1+1)/2$, and 
$\cP(r;2) = \emptyset$ for $r > (d_1+1)/2$.
\end{theorem}
\begin{proof}  
Consider a factor $u$ of $\bs$ with $|u| = r \leq d_1$. As $a_1$
is separating for $\bs$ and $a_1$ occurs in runs of length $d_1$
or $d_1+1$ (inspect the $0$-partition of $\bs$), we have that $u$ is either 
$a_1^r$ or a conjugate of $a_1^{r-1}a_j$ for some $j$, $1< j \leq k$. 
Further, it is evident that there are no squares of words
conjugate to $a_1^{r-1}a_j$, $1< j \leq k$. And, using the same
reasoning for words $u$ of the form $a_1^r$, one determines that
$\bs$ contains the square of $u$ if and only if $2r \leq d_1 + 1$,
in which case there exists exactly one square of each such factor
of length $r$ of $\bs$; namely $(a_1^r)^2$.
\end{proof}

Let $w$ be a factor of $\bs$ with $|w| \in \bigcup_{i=1}^{k} \cD_i(n)$ for some $n$. Roughly speaking, the next theorem shows that if $w^2$ is a factor of $\bs$, then $w$ is a conjugate of a finite product of blocks from the set $\{s_n,s_{n+1},\ldots, s_{n+1-k}\}$, depending on $|w|$ and $d_{n+1}$. For example, if $|w| = r|s_n|$ for some $r$ with $1 \leq r < 1+d_{n+1}/2$, then $w^2 \prec \bs$ if and only if $w$ is one of the first $|s_n|$ conjugates of $s_n^r$. 

\begin{theorem} \label{T:22.11.03(6)}
Let $n$, $r \in \NN^+$.
\begin{itemize}
\item[\emph{(i)}] For $1 \leq r \leq d_{n+1}$,
\begin{equation} \label{eq:power1}
 p(|s_n^r|; 2) = \begin{cases}
                       |s_n| 
                       &\mbox{if $1 \leq r < 1 + d_{n+1}/2$},\\
                       |D_{n-k}| + 1 &\mbox{if $d_{n+1}$ is
                        even and $r = 1 + d_{n+1}/2$}, \\
                        0
                        &\mbox{if $1 + d_{n+1}/2 < r \leq d_{n+1}$}.
                        \end{cases}
\end{equation}
In particular,  
\begin{equation} \label{eq:power2}
 \cP(|s_n^r|; 2) = \begin{cases}
                   \{C_j(s_n^r) ~: ~0 \leq j \leq |s_n| - 1\}
                    &\mbox{if $1 \leq r < 1 + d_{n+1}/2$},\\
                   \{C_j(s_n^r) ~: ~0 \leq j \leq |D_{n-k}|\} 
                    &\mbox{if $d_{n+1}$ is
                        even and $r = 1 + d_{n+1}/2$}, \\
                    \emptyset &\mbox{if $1 + d_{n+1}/2 < r \leq d_{n+1}$}.
                   \end{cases}
\end{equation}
\item[\emph{(ii)}] For $1 \leq r \leq d_{n+1}$ and $i \in [2,k]$ $($with $r\ne d_{n+1}$ if $i=k)$, we have  
\begin{equation} \label{eq:power4}
p(|s_n^rs_{n-1}^{d_n}\cdots s_{n+2-i}^{d_{n+3-i}}s_{n+1-i}|; 2) =
|D_{n+1-i}| + 1.
\end{equation}
In particular, 
\begin{equation} \label{eq:power5}
\cP(|s_n^rs_{n-1}^{d_n}\cdots s_{n+2-i}^{d_{n+3-i}}s_{n+1-i}|; 2)
    = \{C_j(s_n^rs_{n-1}^{d_{n}}\cdots
s_{n+2-i}^{d_{n+3-i}}s_{n+1-i}) ~: ~0 \leq j \leq |D_{n+1-i}|\}.
\end{equation}
\end{itemize}
\end{theorem}

\begin{remark}
For standard Sturmian words $c_\alpha$, we have $s_n = D_{n-1}xy$, where $x$, $y \in
\{a,b\}$ ($x\ne y$), and hence $|D_{n-1}| = q_{n} -2$ for all $n
\geq 1$. Accordingly, Theorem \ref{T:22.11.03(6)} agrees with
Theorem 3 in \cite{dDdL03powe} for the case of a 2-letter alphabet.
\end{remark}

The proof of Theorem \ref{T:22.11.03(6)} requires the following three lemmas.

\begin{lemma} \label{L:09.08.04(1)}  
Let $n \in \NN^+$ and set $u_{i} := s_n^rs_{n-1}^{d_n}\cdots
s_{n+2-i}^{d_{n+3-i}}
     s_{n+1-i}$ for each $i \in [2,k]$ and $1 \leq r \leq d_{n+1}-1$.
Then, for all $i \in [2,k]$, we have \vspace{0.2cm}
\begin{equation} \label{eq:09.08-1}
 \hspace{-6cm} (s_n^{d_{n+1}}s_{n-1}^{d_n}\cdots
    s_{n+2-i}^{d_{n+3-i}}s_{n+1-i})^2D_{n+1-i} \prec \bs,
\end{equation}
and
\begin{equation} \label{eq:09.08-2}
  \hspace{-6cm} u_i^2D_{n+1-i} \prec (s_n^{d_{n+1}}s_{n-1}^{d_n}\cdots
    s_{n+2-i}^{d_{n+3-i}}s_{n+1-i})^2.
\end{equation}
\end{lemma}
\begin{proof}  
Let us first note that, for $i = k$,
$(s_{n}^{d_{n+1}}s_{n-1}^{d_n}\cdots
    s_{n+2-k}^{d_{n+3-k}}s_{n+1-k})^2D_{n+1-k} = s_{n+1}^2D_{n+1-k}$ is a
factor of $\bs$ (by Lemma \ref{L:index1}). Now, for $i \in
[2,k-1]$, by inspection of the $n$-partition of $\bs$, the word
\begin{align*}
 &~ s_n^{d_{n+1}}s_{n-1}^{d_n}\cdots s_{n+2-i}^{d_{n+3-i}}
 s_{n+1-i}s_{n}^{d_{n+1}}s_{n-1}^{d_n}\cdots s_{n+2-k}^{d_{n+3-k}}s_{n+1-k}
 s_n^{d_{n+1}} \\
 = &~ (s_n^{d_{n+1}}s_{n-1}^{d_n}\cdots s_{n+2-i}^{d_{n+3-i}}
 s_{n+1-i})^2s_{n+1-i}^{d_{n+2-i}-1}s_{n-i}^{d_{n-i+1}}\cdots s_{n+1-k}
 s_{n}^{d_{n+1}} \\
 = &~(s_n^{d_{n+1}}s_{n-1}^{d_n}\cdots s_{n+2-i}^{d_{n+3-i}}
 s_{n+1-i})^2D_{n+1-i}D_{n+1-k}^{-1}s_{n}^{d_{n+1}}
\end{align*}
is a factor of $\bs$ (where $D_{n+1-k}$ is a prefix of $s_n$ for
$n\geq k-1$, and $D_{n+1-k}^{-1} = a_{n+2}$ for $1\leq n \leq
k-2$). Thus, assertion \eqref{eq:09.08-1} is proved.

As for the second assertion \eqref{eq:09.08-2}, one can write
\begin{align*}
  &~(s_n^{d_{n+1}}s_{n-1}^{d_n}\cdots s_{n+2-i}^{d_{n+3-i}}s_{n+1-i})^2 \\
 = &~s_n^{d_{n+1}-r}u_i s_n^rs_n^{d_{n+1}-r}s_{n-1}^{d_n}\cdots
      s_{n+2-i}^{d_{n+3-i}}s_{n+1-i} \\
 = &~s_n^{d_{n+1}-r}u_i^2 s_{n+1-i}^{d_{n+2-i}-1}s_{n-i}^{d_{n+1-i}}
    \cdots s_{n-k}s_n^{d_{n+1}-r-1}s_{n-1}^{d_n}\cdots
     s_{n+2-i}^{d_{n+3-i}}s_{n+1-i} \\
 = &~s_n^{d_{n+1}-r}u_{i}^2D_{n+1-i}D_{n-k}^{-1}s_{n}^{d_{n+1}-r-1}
      s_{n-1}^{d_n}\cdots s_{n+2-i}^{d_{n+3-i}}s_{n+1-i},
\end{align*}
which yields the result since $D_{n-k}$ is a prefix of $s_n$ and
$s_{n-1}$ for $n\geq k$, and $D_{n-k}^{-1} = a_{n+1}$ for $1 \leq
n \leq k-1$.
\end{proof}

\begin{lemma} \label{L:02.08.04(2)}  
Let $n \in \NN^+$ and let $u^2 = u^{(1)}u^{(2)}$ be an occurrence
of $u^2$ in $\bs$, where $|s_n| \leq |u| < |s_{n+1}|$.
\begin{itemize}
\item[\emph{(i)}] For all $n\geq 1$, if $|u| = |s_n^r|$ with
$1\leq r \leq d_{n+1}$, then $u^{(1)}$ begins in an $s_n$ block of
the $n$-partition of $\bs$. Moreover, $u^2$ is a factor
of $s_n^{d_{n+1}+2}s_nv^{-1}= s_n^{d_{n+1}+2}D_{n-k}$, where $|v|
= |s_n| - |D_{n-k}|$.
\item[\emph{(ii)}] Let $i \in [2, k-1]$. For all $n \geq i-1$,
if $|u| = |s_n^rs_{n-1}^{d_{n}}\cdots s_{n+2-i}^{d_{n+3-i}}
s_{n+1-i}|$ with $1\leq r \leq d_{n+1}$, then $u^{(1)}$ starts in
an $s_n$ block and contains an $s_{n+1-i}$ block that is followed
by an $s_n$ block in the $n$-partition of $\bs$. Moreover, $u^2$
is a factor of $(s_n^rs_{n-1}^{d_n}\cdots
s_{n+2-i}^{d_{n+3-i}}s_{n+1-i})^2D_{n+1-i}$, which is a
factor of
\[
(s_n^{d_{n+1}}s_{n-1}^{d_n}\cdots
s_{n+2-i}^{d_{n+3-i}}s_{n+1-i})^2 D_{n+1-i}.
\]
\item[\emph{(iii)}] For all $n\geq k-1$, if 
$|u| = |s_n^rs_{n-1}^{d_{n}}\cdots s_{n+2-k}^{d_{n+3-k}}
s_{n+1-k}|$ with $1\leq r \leq d_{n+1}-1$, then 
$u^{(1)}$ starts in an $s_n$ block and contains an $s_{n+1-k}$
block of the $n$-partition of $\bs$. Moreover, $u^2$ is a factor
of $(s_n^rs_{n-1}^{d_n}\cdots
s_{n+2-k}^{d_{n+3-k}}s_{n+1-k})^2D_{n+1-k}$, which is a factor of
\[
  s_{n+1}^2 = (s_n^{d_{n+1}}s_{n-1}^{d_n}\cdots s_{n+2-k}^{d_{n+3-k}}
              s_{n+1-k})^2.
\]
\end{itemize}
\end{lemma}
\begin{proof}  
(i) By similar arguments to those used in the proof of Theorem
\ref{T:22.11.03(5)}, the first claim is obtained from the fact
that $|u| = r|s_n|$ with $1 \leq r \leq d_{n+1}$, together with
Lemma \ref{L:18.05.04(3)}. For the second claim, one uses the fact 
that an $s_n$ block in which $u^{(1)}$ starts is followed by the word
$s_n^ps_{n-1}^{d_n}\cdots s_{n+2-i}^{d_{n+3-i}}s_{n+1-i}s_n$, for
some $i \in [2,k]$ and $0 \leq p \leq d_{n+1}$. Hence, we have
\begin{align*}
s_{n-1}^{d_n}\cdots
s_{n+2-i}^{d_{n+3-i}}s_{n+1-i}s_nG_{n,i-1}^{-1} &=
s_n(s_{n+1-i}^{d_{n+2-i}-1}s_{n-i}^{d_{n+1-i}}
  \cdots s_{n+1-k}^{d_{n+2-k}}s_{n-k})^{-1}s_{n}G_{n,i-1}^{-1} \\
&= s_n(D_{n+1-i}D_{n-k}^{-1})^{-1}D_{n+1-i} \\
&= s_{n}D_{n-k} \\
&= s_{n}s_{n-1}G_{n-1,k-1}^{-1} \\
&= s_{n-1}s_{n}G_{n,1}^{-1} \qquad \qquad \qquad \quad
   \mbox{(by Proposition \ref{P:12.05.04(1)})}\\
&= s_ns_nv^{-1},
\end{align*}
where $|v| = |s_{n}| - |D_{n-k}|$. Therefore, the assertion 
holds provided $u^{(2)}$ does not contain the word
$s_{n-1}s_{n}(w^{-1}G_{n,1})^{-1}$ for some non-empty proper
prefix $w$ of $G_{n,1}$. Indeed, if $s_{n-1}s_{n}(w^{-1}G_{n,1})^{-1} \prec 
u^{(2)}$, then  
$s_{n-1}D_{n-1}w = D_{n-k}\rev{G}_{n,1}D_{n-1}w$ is a factor of 
$u^{(2)}$, where $w \pref
G_{n,1}$.
But this situation is absurd (by Proposition \ref{P:04.08.04(1)}) since this 
word contains a singular $n$-word of the first kind.

(ii) From Lemma \ref{L:18.05.04(3)}, we can argue (as in the proof
of Theorem \ref{T:22.11.03(5)}) that $u^{(1)}$ begins in an
$s_{n+1-i}$ block that is followed by an $s_{n}$ block in the
$n$-partition, or contains an $s_{n+1-i}$ block that is followed
by an $s_n$ block. However, in the first case, we see that $u$
would contain a singular word from the set $\Omega_{n+2-i}^1$,
since $u$ would contain $s_{n+1-i}s_n$, which has $s_{n+1-i}s_{n+2-i}$ as a prefix (see Proposition \ref{P:19.05.04(1)}). This contradicts Proposition \ref{P:04.08.04(1)}, 
and so only the second case can occur. By reasoning as above and
using the fact that $|u| = |s_n^rs_{n-1}^{d_n}\cdots
s_{n+2-i}^{d_{n+3-i}} s_{n+1-i}|$, $u^{(1)}$ must start in the
left-most $s_n$ block in the word
\[
  s_n^rs_{n-1}^{d_n}\cdots s_{n+2-i}^{d_{n+3-i}}s_{n+1-i}s_n,
\]
which appears in the $n$-partition of $\bs$.  
Since $D_{n+1-i} \pref s_{n+2-i} \pref s_n$, $u^{(1)}$ ends within
the first $|D_{n+1-i}|$ letters of the $s_n$ block to the right of
the $s_{n+1-i}$ block. Otherwise, $u$ would contain a singular
word $w \in \Omega_{n+2-i}^1$ which contradicts Proposition  
\ref{P:04.08.04(1)}. To the left of the $s_{n+1-i}$ block, there
is the word $s_n^{d_{n+1}}s_{n-1}^{d_{n}}\cdots
s_{n+2-k}^{d_{n+3-k}}s_{n+1-k}$ and, in view of the fact that $|u|
= |s_{n}^rs_{n-1}^{d_{n}}\cdots s_{n+2-i}^{d_{n+3-i}} s_{n+1-i}|$
with  $1 \leq r \leq d_{n+1}$, one deduces that there exists a $p
\in \NN$ with $0 \leq p \leq |D_{n+1-i}|$  
such that $u^{(1)}$ starts at position $p$ in
$s_{n}^rs_{n-1}^{d_{n}}\cdots s_{n+2-i}^{d_{n+3-i}}s_{n+1-i}s_n$.
This implies $u^2 \prec v^2D_{n+1-i}$  where $v :=
s_{n}^rs_{n-1}^{d_{n}}\cdots s_{n+2-i}^{d_{n+3-i}}s_{n+1-i}$. It
remains to show that $v^2D_{n+1-i}$ is contained in
\[
(s_n^{d_{n+1}}s_{n-1}^{d_n}\cdots
s_{n+2-i}^{d_{n+3-i}}s_{n+1-i})^2 D_{n+1-i} 
\]
which, in turn, is a factor of $\bs$. Indeed, this fact is easily
deduced from Lemma \ref{L:09.08.04(1)}.

(iii) The proof of this assertion is similar to that of (ii), but
with $i = k$ and $1 \leq r \leq d_{n+1}-1$. The details are left to 
the reader.  
\end{proof}


\begin{lemma} \label{L:10.08.04(2)}
For all $n$, $r \in \NN^+$ and $i \in [2, k]$,  
the word $v := s_n^rs_{n-1}^{d_n}\cdots
s_{n+2-i}^{d_{n+3-i}}s_{n+1-i}$ is primitive. 
\end{lemma}
\begin{proof} Suppose on the contrary that the given word $v$ is not primitive, i.e., suppose $v = u^p$
for some non-empty word $u$ and integer $p \geq 2$. Then 
$|v|_{a_j} = p|u|_{a_j}$ for each letter $a_j \in \cA_k$, i.e., $p$ divides 
$|v|_{a_j}$ for each $a_j \in \cA_k$. In particular, $p$
divides
\begin{align*}
  |v|_{a_1} &= r(Q_n - P_n) + d_{n}(Q_{n-1} - P_{n-1}) + \cdots + (Q_{n+1-i} - P_{n+1-i}) \\ 
  &= |v| - (rP_n + d_{n}P_{n-1} + \cdots + d_{n+3-i}P_{n+2-i} + P_{n+1-i}),
\end{align*}
by Proposition \ref{P:13.09.04(1)}. Thus, $p$ must also divide
$rP_n + d_{n}P_{n-1} + \cdots + d_{n+3-i}P_{n+2-i} + P_{n+1-i}$.
But gcd$(|v|, rP_n + d_{n}P_{n-1} + \cdots + d_{n+3-i}P_{n+2-i} +
P_{n+1-i}) = 1$, which yields a contradiction; whence $p=1$, and
therefore $v$ is primitive. 
\end{proof}


\begin{proof}[Proof of Theorem $\ref{T:22.11.03(6)}$]  
We simply prove \eqref{eq:power2} and  \eqref{eq:power5}  
as the elements of the respective sets are mutually distinct (by Lemma
\ref{L:10.08.04(2)}), which implies \eqref{eq:power1} and 
\eqref{eq:power4}.  

(i) As we have already seen, for each $i \in [2,k]$,  
$v := s_{n}^{d_{n+1}}s_{n-1}^{d_n}\cdots
s_{n+2-i}^{d_{n+3-i}}s_{n+1-i}s_nG_{n,i-1}^{-1} = s_n^{d_{n+1}+2}D_{n-k}$
is a factor of $\bs$. 
Thus, by Lemma \ref{L:02.08.04(2)}(i), it suffices to
find all squares of words $u$ with $|u| = r|s_n|$ ($1\leq r
\leq d_{n+1}$) that occur in the word $v$. In fact, one need only
consider occurrences of $u^2$ starting in the left-most $s_n$ block of
$v$, and the result easily follows.

(ii) Let $i \in [2, k-1]$. By Lemma \ref{L:09.08.04(1)}, the word
$(s_{n}^{d_{n+1}}s_{n-1}^{d_n}\cdots
s_{n+2-i}^{d_{n+3-i}}s_{n+1-i})^2$ is a factor of $\bs$ and it
contains the word
\[
  v: = (s_{n}^{r}s_{n-1}^{d_n}\cdots s_{n+2-i}^{d_{n+3-i}}s_{n+1-i})^2
       D_{n+1-i},
\]
for any $r$ with $1 \leq r \leq d_{n+1}$.  
Therefore, $u^2$ is a factor of $\bs$ for each word $u$ given by
\[
  u := C_j(s_{n}^{r}s_{n-1}^{d_n}\cdots s_{n+2-i}^{d_{n+3-i}}s_{n+1-i}),
\quad 0\leq j \leq |D_{n+1-i}|.
\]
Conversely, if $u^2 \prec \bs$ with $|u| =
|s_{n}^{r}s_{n-1}^{d_n}\cdots s_{n+2-i}^{d_{n+3-i}}s_{n+1-i}|$ for
some $1 \leq r \leq d_{n+1}$, then $u^2 \prec v$, by Lemma
\ref{L:02.08.04(2)}(ii). And, since $|v| = 2|u| + |D_{n+1-i}|$, we must have 
$u = C_j(s_n^rs_{n-1}^{d_{n}}\cdots
s_{n+2-i}^{d_{n+3-i}} s_{n+1-i})$ for some $j$ with $0\leq j \leq
|D_{n+1-i}|$.

The case $i = k$ is proved similarly, using Lemma \ref{L:02.08.04(2)}(iii).
\end{proof}


\subsection{Cubes and higher powers} \label{SS:higher}

Our subsequent analysis of cubes and higher powers occurring in $\bs$
is now an easy task due to the above consideration of squares. 
Extending Theorem \ref{T:22.11.03(6)} (see
Theorem \ref{T:11.08.04(2)} below), only requires the following
lemma, together with arguments used in the proof of Theorem
\ref{T:22.11.03(6)}.

\begin{lemma} \label{L:11.08.04(1)} 
Let $n \in \NN^+$ and suppose $u^3 \prec \bs$ with $|s_n| \leq |u| < |s_{n+1}|$. Then
$u^3$ does not contain a singular word from the set
$\Omega_{n+1-i}^1$ for any $i \in [1,k-1]$.
\end{lemma}
\begin{proof}  
Suppose on the contrary that $u^3 = u^{(1)}u^{(2)}u^{(3)}$
contains a singular word $w \in \Omega_{n+1-i}^1$ for some $i \in
[1,k-1]$. By Proposition \ref{P:04.08.04(1)}, $w$ is not a factor of
$u^{(3)}$, and therefore every occurrence of $w$ must begin in
$u^{(1)}$ or $u^{(2)}$, both of which are followed by $u$ again. 
Accordingly, there exists a $p \in \NN$
such that $w$ starts at position $p$ in both $u^{(1)}$ and
$u^{(2)}$. Reasoning, as in the proofs of Propositions  
\ref{P:27.11.03(2)} and \ref{P:04.08.04(1)}, yields the
contradiction $|u| \geq |s_{n+1}|$.
\end{proof}

\begin{theorem} \label{T:11.08.04(2)}
Let $n$, $r$, $l \in \NN^+$, $l \geq 3$.
\begin{itemize}
\item[\emph{(i)}] For $1 \leq r \leq d_{n+1}$,
\begin{equation} \label{eq:power11}
 p(|s_n^r|; l) = \begin{cases}
                       |s_n| 
                       &\mbox{if $1 \leq r < (d_{n+1}+2)/l$},\\
                       |D_{n-k}| + 1 &\mbox{if
                                            $r = (d_{n+1}+2)/l$}, \\
                        0
                        &\mbox{if $(d_{n+1}+2)/l < r \leq d_{n+1}$}.
                        \end{cases}
\end{equation}
In particular,  
\begin{equation} \label{eq:power22}
 \cP(|s_n^r|; l) = \begin{cases}
                       \{C_j(s_n^r) ~: ~0 \leq j \leq |s_n| - 1\}
                       &\mbox{if $1 \leq r < (d_{n+1}+2)/l$},\\
                       \{C_j(s_n^r) ~: ~0 \leq j \leq |D_{n-k}|\}
                       &\mbox{if $r = (d_{n+1}+2)/l$}, \\
                       \emptyset
                       &\mbox{if $(d_{n+1}+2)/l < r \leq d_{n+1}$}.
                        \end{cases}
\end{equation}
\item[\emph{(ii)}] For $1\leq r \leq d_{n+1}$ and $i \in [2, k]$
$($with $r \ne d_{n+1}$ if $i = k)$, we have
\begin{equation} \label{eq:power44}
p(|s_n^rs_{n-1}^{d_n}\cdots s_{n+2-i}^{d_{n+3-i}}s_{n+1-i}|; l) =
0.
\end{equation}
\end{itemize}
\end{theorem}
\begin{proof} 
(i) Suppose $u \prec \bs$ with $|u| = r|s_n|$ for some $r$, $1\leq r\leq d_{n+1}$, 
and consider an occurrence of $u^{l} =
u^{(1)}u^{(2)}\cdots u^{(l)}$ in $\bs$, $l \geq 3$. By Lemma
\ref{L:02.08.04(2)}(i), $u^{(1)}u^{(2)}$ begins in an $s_n$ block
of the $n$-partition of $\bs$, and by Lemma \ref{L:11.08.04(1)},
$u^l$ does not contain a singular $n$-word of the first kind. So,
as in the proof of Theorem \ref{T:22.11.03(6)}(i), one infers that
$u^l$ is contained in the word $v := s_{n}^{d_{n+1}+2}D_{n-k}
\prec \bs$, and the rest of the proof follows in much the same fashion. 

(ii) In this case, assume $u^3 = u^{(1)}u^{(2)}u^{(3)}$ occurs in
$\bs$. By (ii) and (iii) of Lemma \ref{L:02.08.04(2)}, $u^{(1)}$ begins in an $s_n$ block and contains an
$s_{n+1-i}$ block that is followed by an $s_n$ block in the
$n$-partition of $\bs$. Accordingly, $u^3$ contains the word
$s_{n+1-i}s_n$, and hence contains a singular word $w \in
\Omega_{n+2-i}^1$, which contradicts Lemma \ref{L:11.08.04(1)}.
\end{proof}

\subsection{Examples}

\begin{example} \label{Ex:13.08.04}
Let us demonstrate Theorems \ref{T:22.11.03(6)} and
\ref{T:11.08.04(2)} with an explicit example. Consider the
standard episturmian word $\bs$ over $\cA_3 = \{a_1,a_2,a_3\}
\equiv \{a,b,c\}$ with periodic directive word $(abcca)^\omega$.
We have $(d_{n})_{n\geq1} = (1,1,2,\ov{2,1,2})$, and hence
\begin{align*}
s_1 &= ab, \quad |s_1| = 2;  \\
s_2 &= abac, \quad |s_2| = 4; \\
s_3 &= abacabacaba, \quad |s_3| = 11; \\
s_4 &= abacabacabaabacabacabaabacabacab,\quad |s_4| = 32; \\
s_5 &= abacabacabaabacabacabaabacabacababacabacabaabacabacabaabac,
       \quad |s_5| = 58.
\end{align*}
Also, $D_0 = \empt$, $D_1 = a$, $D_2 = abacaba$, and $D_3 =
abacabacabaabacabacaba$.

We shall simply consider squares and cubes of words of length $m$
occurring in $\bs$ with $|s_3| \leq m < |s_6|$. By Theorem
\ref{T:22.11.03(5)}, we need only consider lengths $m$ in the  
set
\[
  \cS := \{|s_3|, 2|s_3|, |s_4|, |s_5|, 2|s_5|, |s_3s_2|,
                     |s_3^2s_2|, |s_3s_2^2s_1|,|s_4s_3|, |s_5s_4|, |s_5^2s_4|,
                     |s_5s_4s_3|\}.
\]

According to Theorem \ref{T:22.11.03(6)}(i), for $3 \leq n \leq
5$, we have
\begin{align*}
&p(|s_3|;2) = |s_3| = 11, \\
&p(2|s_3|;2) = |D_0| + 1 = 1, \\
&p(|s_4|;2) = |s_4| = 32, \\
&p(|s_5|;2) = |s_5| = 58, \\
&p(2|s_5|;2) = |D_2| + 1 = 8.
\end{align*}
Also, part (ii) of Theorem \ref{T:22.11.03(6)} gives 
\begin{align*}
&p(|s_3s_2|;2) = |D_2| + 1 = 8, \\
&p(|s_3^2s_2|;2) = |D_2| + 1 = 8, \\
&p(|s_3s_2^2s_1|;2) = |D_1| + 1 = 2,\\
&p(|s_4s_3|;2) = |D_3| + 1 = 23, \\
&p(|s_5s_4|;2) = |D_4| + 1 = 34, \\
&p(|s_5^2s_4|;2) = |D_4| + 1 = 34, \\
&p(|s_5s_4s_3|;2) = |D_3| + 1 = 23.
\end{align*}
Furthermore, from Theorem \ref{T:11.08.04(2)}, one has
\begin{align*}
  &p(|s_3|;3) = |s_3| = 11, \\
  &p(2|s_3|;3) = 0, \\
  &p(|s_4|;3) = |D_1| + 1 = 2, \\
  &p(|s_5|;3) = |s_5| = 58, \\
  &p(2|s_5|;3) = 0,
\end{align*}
and $p(m;3) = 0$ for all other lengths $m \in \cS$.

For instance, the sole factor of $\bs$ of length $2|s_3| = 22$
that has a square in $\bs$ is
\[
  s_3^2 = abacabacabaabacabacaba,
\]
and the eight squares of length $2|s_3s_2| = 30$ are the squares
of the first eight conjugates of $s_3s_2 = abacabacabaabac$;
namely
\[
  (abacabacabaabac)^2, ~~(bacabacabaabaca)^2,  
   ~~\ldots~~, ~~(cabaabacabacaba)^2.
\]
The only factors of length $|s_4|=32$ that have a cube in $\bs$ 
are the first two conjugates of $s_4$, i.e.,
\[
  s_4^3 \prec \bs \quad \mbox{and} \quad (C_1(s_4))^3 = (a^{-1}s_4a)^3
\prec \bs.
\]
\qed
\end{example}

\begin{example} 
The \emph{$k$-bonacci word} is the standard episturmian word 
$\boldsymbol{\eta}_k \in \cA_k^\omega$ with directive word 
$(a_1a_2\cdots a_k)^\omega$. 
Since all $d_i = 1$, we have $s_n = s_{n-1}s_{n-2}\cdots s_{n-k}$ for 
all $n \geq 1$ (and the lengths $|s_n|$ are the \emph{$k$-bonacci numbers}). 
Thus, for fixed $n \in \NN^+$ and $l\geq2$, if $w^l \prec \boldsymbol{\eta}_k$ 
with $|s_n| \leq |w| < |s_{n+1}|$, then we necessarily have 
$|w| = |s_{n}| + |s_{n-1}| + \cdots + |s_{n+1-i}|$ 
for some $i \in [1,k-1]$ (by Theorem \ref{T:22.11.03(5)}). The preceding 
main theorems reveal that
\[
  \cP(1;2) = \{a_1\}, \quad \cP(|s_n|;2) = \cC(s_n) = \Omega_n^0 \quad 
  \mbox{and} 
  \quad \cP(|s_n|;3) = \{C_j(s_n) ~: ~0\leq j \leq |D_{n-k}|\}.
\]
Furthermore, for each $i \in [2,k-1]$, we have 
\[
  \cP(|s_ns_{n-1}\cdots s_{n+1-i}|;2) 
  =\{C_j(s_ns_{n-1}\cdots s_{n+1-i}) ~: ~0\leq j \leq |D_{n+1-i}|\}.
\]
All other $\cP(|w|;l) = \emptyset$, $l\geq2$. In particular, 
$k$-bonacci words are $4$-power free.
\qed
\end{example}

\section{Concluding remarks}

Using the results of Section 6, it is possible to determine the exact number of distinct squares in each  building block $s_n$, which extends Fraenkel and Simpson's result \cite{aFjS99thee} concerning squares in the \emph{finite Fibonacci words}. Such work forms part of the present author's PhD thesis \cite[Chapters 6 and 7]{aG06onst}.  

Theorems \ref{T:02.08.04(1)}, \ref{T:22.11.03(6)} and
\ref{T:11.08.04(2)} also suffice to describe all integer powers occurring in 
any (episturmian) word $\bt \in \cA_k^\omega$ that is equivalent to 
$\bs$. (See \cite[Theorem 3.10]{jJgP02epis} for a definition of such $\bt$.)  
The problem of determining all integer powers occurring in general 
standard episturmian words (with not all $d_i$ necessarily positive) 
remains open.

\section{Acknowledgements}

The author would like to thank the anonymous referee for their helpful comments and suggestions.


\end{document}